\documentclass[12pt,a4paper]{article}

\usepackage{amsmath}
\usepackage{comment}
\usepackage{amsfonts}
\usepackage{graphicx}
\usepackage{color}
\usepackage{tcolorbox}
\usepackage{hyperref}
\graphicspath{{./figures/}{/figures-new/}}
\usepackage{enumitem}
\usepackage{relsize}
\usepackage{tikz}
\usetikzlibrary{arrows.meta}
\usepackage{xcolor}
\definecolor{pathgreen}{RGB}{0,190,65}
\definecolor{labelblue}{RGB}{0,0,255}
\definecolor{pointred}{RGB}{255,0,0}
\definecolor{axisgray}{RGB}{145,145,145}
\definecolor{boxgray}{RGB}{238,238,238}%

\DeclareMathOperator\sign{sign}

\newcommand{\C}{\ensuremath{\mathbb{C}}}

\newcommand{\rme}{\mathrm{e}}
\newcommand{\rmi}{\mathrm{i}}

\renewcommand{\Im}{\operatorname{Im}}
\renewcommand{\Re}{\operatorname{Re}}

\newtheorem{prop}{Proposition}
\newtheorem{defn}{Definition}
\newtheorem{rem}{Remark}

\title{Singularly embedded solitons\\ and Stokes phenomenon\\ in the fifth order KdV equation}
\author{R.~Barrio$^{1}$, 
        V.~Gelfreich$^{3}$\footnote{Corresponding author e-mail: \tt v.gelfreich@warwick.ac.uk}, J.T.~Lázaro$^{4}$, J.R.Pacha$^{4}$\\[1,0em]
\small $^{1}$Departamento de Matemática Aplicada, University of
       Zaragoza\\
\small  C.~Pedro Cerbuna 12, E-50009 Zaragoza, Spain\\ [0,5em]
\small $^{3}$Mathematics Institute, University of Warwick\\
\small Coventry CV4 7AL, United Kingdom\\ [0,5em]
\small $^{4}$Departament de Matemàtiques, Universitat
Politècnica de Catalunya\\
\small Av.~Diagonal 647, 08028 Barcelona, Spain
}

\date{November 12, 2026}

\begin{document}
\maketitle

\begin{abstract}
We study embedded solitons (ES) in a fifth-order Korteweg–de Vries (KdV) equation  near a reversible Hamiltonian saddle–centre bifurcation. The solitons correspond to homoclinic orbits to an equilibrium whose linearisation has a pair of real and a pair of purely imaginary eigenvalues. Existence of symmetric
homoclinic near the bifurcation is a codimension-one phenomenon governed by a parameter-dependent Stokes constant associated with an inner equation.  Using complex matching, we derive an oscillatory large-parameter asymptotic  formula for the Stokes constant. This formula implies the existence of a countable sequence of zeros and gives their asymptotic distribution in parameter space. We also develop two complementary high-precision numerical procedures for evaluating the Stokes constant: one based on the late terms of the formal inner expansion and another on direct computation of the stable and unstable manifolds. The numerical results agree with each other and with the asymptotic prediction, confirming both the existence and the predicted locations of the embedded-soliton branches.
\end{abstract}

\tableofcontents

\newpage

\section{Introduction\label{sec:intro}}

\subsection{Embedded solitons in a fifth order KdV}

A soliton is a solution of a non-linear partial differential equation (PDE) that 
is localised in space and moves with constant speed while preserving its shape.
In some equations the soliton may coexist with linear harmonic waves that travel with
the same velocity \cite{ChampneysGroves97}. 
Those solutions are  called embedded solitons (ES)
to stress their relation with the  spectrum of the linearised equation \cite{YangEmbedded}. 
ES are observed in  various models, for example, in models of non-linear optics, capillary-gravity surface wave and higher-order KdV type equations~\cite{Champmeys01PhysD}.

It is known that  emergence of ES can be related to a bifurcation in
the PDE's dispersion relation that establishes connection between
wave number $k$ and phase velocity $c$ of a small-amplitude harmonic solution to the  problem
linearised around an equilibrium state. Following \cite{ChampneysGroves97}
 we illustrate this situation  with the help of  the 5th order KdV model 
\begin{equation}\label{Eq:PDE}
u_t=\Bigl[\tfrac2{15}u_{xxxx}-bu_{xx}+au+\tfrac32u^2 -\tfrac12u_x^2+\left(uu_x\right)_x
\Bigr]_x.
\end{equation}
This equation was derived from the Euler equation for water waves with surface tension 
using the regular Hamiltonian perturbation theory~\cite{CraigGroves94}.
The equation depends on two real parameters $a$ and $b$.
Looking for a solution in the form $u=A\mathrm{e}^{ik(x-ct)}$
and neglecting the non-linear terms we obtain  the dispersion relation
in the form
\begin{equation}\label{Eq:DispRel}
\tfrac2{15}k^4+bk^2+a+c=0\,.
\end{equation}
This equation determines the phase speed $c$ of small amplitude 
harmonic waves as a function of the wave number $k$. The dispersion relation can also be considered as an equation for possible values of $k$ for a fixed $c$. 
Real values of $k$ correspond to small amplitude
harmonic waves, while non-real ones correspond to 
exponentially increasing/decreasing solutions of the linearised equation and open a possibility for existence of solitons. 
It is easy to see that for $a+c<0$  the equation \eqref{Eq:DispRel}
has a pair of real
solutions and a pair of purely imaginary ones. In this region
of parameters a soliton
may co-exist with a small amplitude harmonic wave. 

A bifurcation happens at $a+c=0$ where a long linear wave with  $k\ll1$, which exists for $c>-a$,
becomes an unstable linear mode for $c<-a$.
The numerical study \cite{ChampneysGroves97} suggests that the values of parameters
corresponding to ES form a system of isolated lines 
on the $(b,a+c)$-plane which emerge from
a sequence of points with $b=b_k<0$ and $a+c=0$. 

Looking for a travelling wave solution of the equation \eqref{Eq:PDE} 
in the form $u=u(x-ct)$,  integrating once and setting 
the integration constant to zero, we obtain an ordinary differential equation (ODE)
\begin{equation}
\label{Eq:ODE4wave}
\tfrac2{15}u^{(iv)}-bu''+ c u+\tfrac32u^2 -\tfrac12(u')^2+\left(uu'\right)'=0
\end{equation}
where the primes refer to derivatives with respect to the new independent variable $\tau=x-ct$. 
We also substituted $a=0$ to shorten the notation without loosing in generality as the equation depends only on the sum $a+c$.
Solitons correspond to 
solutions of this equation that converge to zero at infinity. 

The equation~\eqref{Eq:ODE4wave} inherits reversibility from the original PDE: if a function $u(\tau)$ is a solution, then so is 
$u(-\tau)$. Moreover, \eqref{Eq:ODE4wave} can be rewritten as a Hamiltonian system with two degrees of freedom, where the trivial solution $u \equiv 0$ corresponds to an equilibrium point.
This reformulation of the problem is described in Appendix~\ref{Se:hamiltonian}. 

In dynamical systems theory, a solution of an ODE that converges to the same equilibrium as $\tau \to \pm\infty$ is called a homoclinic orbit. Consequently, homoclinic solutions to the origin correspond to solitons.
    The theory of dynamical systems provides a range of tools suitable for studying homoclinic solutions, including regular and singular perturbation theories, bifurcation theory, as well as various numerical methods (see \cite{ChampneysMalomedYang} for a detailed discussion, as well as \cite{Meyer,Kuznetsov}).

\subsection{Bifurcation diagram}

\begin{figure}
	\centerline{\includegraphics[scale=1.6]{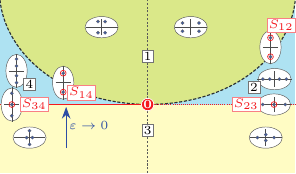}}
	\caption{\protect{\small Classification of the different regions
			of the linearisation at the origin of the 
			equation~\eqref{Eq:ODE4wave} on the $(b,c)$-plane.
			\label{classification}}}
\end{figure}

The linear stability of the equilibrium is determined by the characteristic exponents $\lambda = \mathrm i k$, where $k$ satisfies the dispersion relation~\eqref{Eq:DispRel}, therefore
\[
\frac{2}{15}\lambda^4 - b \lambda^2 + c = 0.
\]
This is a biquadratic equation and, in general,  admits four distinct solutions. As the parameters vary, two of these solutions may coincide. A straightforward calculation shows that this occurs when 
the determinant $b^2 - \frac{8}{15}c = 0$ or when $\lambda^2=0$ for $c = 0$. These two lines partition the parameter plane into four regions
sketched in Fig.~\ref{classification}.

The equation~\eqref{Eq:ODE4wave} admits no soliton solutions when the parameters $(b,c)$ lie in Region~4, since the characteristic exponents are purely imaginary and therefore there are no solutions that  decay 
exponentially when  $\tau \to \pm\infty$. In contrast, solitons may exist in all other regions.

From a generic perspective, one expects solitons to occur for an open set of parameter values in Regions~1 and~2. This is because the stable and unstable manifolds are two-dimensional surfaces contained within a three-dimensional level set of the Hamiltonian \cite{ChampneysGroves97,Champmeys01PhysD}.

The situation is different in Region~3, where the stable and unstable manifolds of the equilibrium are one-dimensional and, consequently, the existence of homoclinic solutions is expected to be a codimension-one phenomenon.
The embedded solitons arise for parameter values in Region~3 and correspond to solutions homoclinic to the saddle--centre equilibrium at the origin~\cite{ChampneysGroves97,Champmeys01PhysD}.

\subsection{Existence of solitons}

The above arguments do not, however, establish the existence of solitons. When the parameters lie on the curve
$c = \frac{3}{5}(2b + 1)(b - 2)$ with $ -\frac{1}{2} < b < 2$,
homoclinic solutions can be expressed explicitly in terms of elementary functions~\cite{ChampneysGroves97}. Apart from this special case, no explicit soliton solutions of \eqref{Eq:ODE4wave} are known.
Additional soliton solutions can be obtained near the boundaries between regions using bifurcation theory~\cite{ChampneysMalomedYang}.
The bifurcation theory employs normal forms to describe changes in the dynamics that occur when parameters cross a bifurcation curve.  
In our case, the normal forms are integrable and the resulting dynamics can be analysed explicitly. In particular, the normal forms predict the emergence of a homoclinic solution from the equilibrium when parameters move from Region~3 to Region~2, from Region~4 to Region~1, and from Region~4 to Region~3.

However, the normal forms theory provides only an approximate description of the original dynamics. In the first two cases, regular perturbation theory can be used to show that the existence of a homoclinic solution in the normal form implies the existence of a corresponding homoclinic solution in the full system.

The situation near $S_{34}$ is different, and the analysis of persistence of homoclinic orbits requires singular perturbation theory~\cite{Lombardi00,IoossLombardi04,GelfreichLerman,TovbisPelinovsky06}.
The main conclusion of the theory is that, when the parameters lie on the curve $S_{34}$, the equation~\eqref{Eq:ODE4wave} with $c=0$,
called the inner equation,
\begin{equation}
	\label{Eq:ODE4waveA=0}
	\frac{2}{15}u^{(iv)} - b u'' + \frac{3}{2}u^2 - \frac{1}{2}(u')^2 + (u u')' = 0,
\end{equation}
can be used to define a Stokes constant $S = S(b)$. This constant encodes the splitting of stable and unstable manifolds of the degenerate saddle--centre equilibrium and plays a central role in determining whether homoclinic orbits persist under perturbation.

If $S(b) \neq 0$, then there are no single-bump homoclinic solutions for parameter values near the point $(b,0)$ in Region~3. Consequently, a primary homoclinic solution can only exist near a zero of the Stokes constant. In the context of equation~\eqref{Eq:ODE4wave}, Tovbis and Pelinovsky~\cite{TovbisPelinovsky06} showed that every simple zero $b_0$, $S(b_0)=0$, generates a family of homoclinic orbits, which form a smooth curve in the $(b,c)$-plane and starts at $(b_0,0)$.

%

 The Stokes constant $S$ is an analytic function of the parameter, although no explicit closed-form expression is known and usually its values
 are obtained numerically.

We  show that, for large positive $\gamma=-5b/2$ (or equivalently for large negative $b$), the inner problem \eqref{Eq:ODE4waveA=0}
can itself be interpreted within the framework of singular perturbation theory. In this setting, two complex singularities of an outer problem contribute comparably, leading to oscillatory asymptotics 
\begin{equation}\label{Eq:stokesasympt}
S_0(\gamma)\approx  \mathrm i\,
2\left(\tfrac34\right)^{\frac 14}|K_0|
\gamma^{-{23}/{16}}\sin\left(
\frac{\pi}{16}+4\sqrt{3\gamma}
\right)
\end{equation}
where $K_0$ is a Stokes constant obtained from a parameterless
problem described in Section~\ref{Se:inner4inner},
$K_0=0.251346943975660732\ldots \cdot \mathrm e^{\mathrm i \pi/16}$.
It  follows that, if $\gamma$ is sufficiently large,
 there is a zero of $S_0(\gamma)$ 
near every point 
\begin{equation}\label{Eq:gamma0asympt}
\gamma_0(n)=\left(\frac{\pi n+\frac{15\pi}{16}}{4\sqrt{3}}\right)^2,
\end{equation}
where the approximation for $S_0(\gamma)$ changes its sign.
It is interesting to note that our numerical experiments 
show that every interval $(\gamma_0(n),\gamma_0(n+1))$ 
contains exactly one zero $\gamma_n^*$
for every $n\in \mathbb N$, and 
 that there are no other zeroes.
 Moreover, the difference 
$\sqrt{\gamma_n^*}-\sqrt{\gamma_0(n)}$
monotonically converges to zero.

This is a rather curious example of the theory of asymptotic beyond all orders. Usually the theory is used to estimate exponentially small terms which contain a pre-factor with 
the imaginary part of a leading singularity in the exponent. In our case the leading singularity 
  is on the real axis
and, consequently, the ``exponentially small" pre-factor
equals to one.
 
 To confirm our theoretical predictions, we numerically computed the values of $S_0$ and analysed the locations of its zeros using two complementary methods. The first is based on the analytic continuation of solutions, while the second relies on an analysis of the coefficients of formal solutions. In these numerical experiments, high-precision arithmetic is essential for obtaining reliable results, particularly for larger values of $\gamma$, where numerical evaluation becomes increasingly delicate.

The rest of the paper is organized in the following way.
\begin{itemize}
    \item Section~\ref{Se:bifurcationSingPert}: summary of the singular perturbation theory, inner and outer limits for the equation \eqref{Eq:ODE4wave}, reduction of persistence of solitons     to zeroes of the Stokes constant
    \item Section~\ref{Se:inner}: analysis of the inner problem, definition and basic properties of the Stokes constant
    \item Section~\ref{Se:Sasympt}: inner problem as a singularly perturbed problem, outer and inner equations for the inner equation, asymptotics of $S_0$ for large $\gamma$,
    existence of infinitely many zeroes.
    \item Section~\ref{Se:numerics}: numerical methods for evaluation of the Stokes constant and finding its zeroes.

\end{itemize}

\section{Hamiltonian bifurcation and singular perturbations\label{Se:bifurcationSingPert}}

\subsection{Scaling invariants of the equation}

The study of  solitons emerging at the bifurcation relies
on various changes of scales in the ODE. Since these
substitutions affect the coefficients without changing the form
of the equation, it is convenient to consider a  more general equation
\begin{equation}\label{Eq:mainODE}
	a_2u^{(iv)}+a_1u''+ a_0 u+b_0 u^2+b_1 ({u'}^2+2u u'')=0
\end{equation}
which depends on five real parameters $a_0,a_1,a_2,b_0,b_1$.
We assume that the quadratic part is not degenerate, i.e., $b_0,b_1\ne0$.
Equation~\eqref{Eq:mainODE} is an ODE of order four with a quadratic nonlinearity,
which coincides with \eqref{Eq:ODE4wave} when $a_2=\frac2{15}$,
$a_1=-b$, $a_0=c$, $b_0=\frac32$ and $b_1=\frac12$.
In spite of the substantially larger number of parameters
 equation~\eqref{Eq:mainODE} is only a minor generalisation
of~\eqref{Eq:ODE4wave}.
Indeed, if $a_2,b_0,b_1\ne0$, then 
the substitution $u(t)=\alpha\tilde u(\tau)$ 
 with $t=\beta\tau$,
\[
\alpha=\frac{a_2b_0}{b_1^2}
\quad\text{and}\quad
\beta=\left|\frac{b_0}{b_1}\right|^{1/2}
\]
transforms  \eqref{Eq:mainODE} into a standard form
\begin{equation}\label{Eq:standard form}
u^{(iv)}+p_1 u''+p_2 u+u^2+\sigma_0({u'}^2+2uu'')=0
\end{equation}
where $\sigma_0=\sign(b_0b_1)$,
$
p_1=\frac{a_1|b_1|}{a_2|b_0|}
$
and
$
p_2=\frac{a_0b_1^2}{a_2b_0^2}
$. 
For example, for  equation~\eqref{Eq:ODE4wave} we get
$p_1=-5b/2$, $p_2=5c/6$ and $\sigma_0=1$.
In equation \eqref{Eq:standard form}, we have dropped the tildes for brevity. Here and later, we use primes to denote derivatives, regardless of the letter used for the function’s argument.

We call  parameters $p_1$ and $p_2$ {\em scaling invariants}, since any two equations~\eqref{Eq:mainODE} with non-vanishing coefficients can be transformed into one another if and only if their scaling invariants coincide. These invariants provide a convenient criterion for establishing the equivalence of two equations without explicitly constructing the transformation between them.

\subsection{Hamiltonian bifurcation and singular perturbations}

In the search for embedded solitons, we consider the reversible Hamiltonian bifurcation which happens when the linear stability of the equilibrium changes along the line $S_{34}$ of Fig.~\ref{classification}. An embedded soliton can exists only
in the area below the line, where the equilibrium is of  the saddle--centre type. That means $a_2a_0<0$ and, equivalently, $p_2<0$.
The line $S_{34}$ belongs to the region where $a_2a_1>0$ and, equivalently, $p_1>0$.
We denote
the  eigenvalues of the equilibrium by $\pm\lambda$ and $\pm \rmi\,\omega$  assuming that $\omega,\lambda>0$. 
On the line $S_{34}$, $\lambda=0$ is a double eigenvalue and $\omega>0$.

The scaling of the time variable changes the values of $\omega$ and $\lambda$ without affecting the ratio $\lambda/\omega$. 
This freedom allows us to consider the corresponding equations in two
equivalent ways.

First, 
 an equation~\eqref{Eq:mainODE} with $a_2a_1>0$, $a_2a_0\le 0$ and $b_0b_1\ne0$ has the scaling invariants $p_1>0$ and $p_2\le0$
and, consequently, it can be transformed into an equation of the form
considered in \cite{TovbisPelinovsky06},
\begin{equation}\label{Eq:bifurcform}
u^{(iv)} + (1-\varepsilon^2) u'' - \varepsilon^2 u + (1-\varepsilon^2)u^2 + \gamma \bigl({u'}^2 + 2u u''\bigr) = 0
\end{equation}
with $\varepsilon\in[0,1)$ and $\gamma\ne 0$.
 Indeed, computing the scaling invariants
for this equation we get $p_1=|\gamma|$ and $p_2=-\varepsilon^2\gamma^2/(1-\varepsilon^2)^2$.
Thus we cover all possible values  of the scaling invariants 
in the region where $p_1>0$ and $p_2\le 0$.

In the case of equation~\eqref{Eq:bifurcform},
the characteristic equation 
 takes the form
 $\lambda^4+(1-\varepsilon^2)\lambda^2-\varepsilon^2=0$
 and consequently
   $\lambda=\varepsilon$ and $\omega=1$.  Thus in this formulation, the problem reduces to the study of a reversible Hamiltonian $(0^2,\pm i\omega)$ bifurcation for  equation~\eqref{Eq:bifurcform} which happens at $\varepsilon=0$. 

Second,  
the standard scaling $u(t)=\varepsilon^2 v(\tau)$ with $\tau=\varepsilon t$
transforms  equation~\eqref{Eq:bifurcform} into a singularly perturbed problem,
\begin{equation}\label{Eq:sing}
\varepsilon^2 v^{(iv)} + (1-\varepsilon^2) v'' - v + (1-\varepsilon^2)v^2 + \varepsilon^2 \gamma \bigl({v'}^2 + 2v v''\bigr) = 0,
\end{equation}
where the small parameter multiplies the highest derivative.
In this formulation, we get $\lambda=1$ and $\omega=\varepsilon^{-1}$
goes to $+\infty$ when $\varepsilon \to 0$. 
In the limit $\varepsilon \to 0$  equation~\eqref{Eq:sing} reduces to an equation of the second order,
\[
v'' - v + v^2 = 0.
\]
In the singular perturbation theory, this equation
is called an {\em outer equation}.
The outer equation admits an explicit homoclinic solution
\begin{equation}\label{Eq:v00}
v_0(\tau)=\frac{3}{2}\cosh^{-2}\!\left(\frac{\tau}{2}\right).
\end{equation}
After inverting the scalings, we obtain an approximation for a soliton
of the original equation~\eqref{Eq:ODE4wave} provided it survives 
the influence of the singular perturbation.

\subsection{Persistence of a soliton is of co-dimension one\label{Se:solitons}. }

Despite the existence of the homoclinic orbit 
in equation~\eqref{Eq:sing} with $\varepsilon=0$, such an orbit is generically destroyed when $\varepsilon>0$, see \cite{Lombardi00,GelfreichLerman}. To study the persistence of this solution for $\varepsilon \neq 0$, we consider the stable and unstable solutions of equation~\eqref{Eq:bifurcform}.
Since the unstable direction at the equilibrium is one dimensional, 
 the unstable solution, $u_-(t)\to 0$ as $t\to -\infty$, 
 is defined uniquely up to 
a constant shift in time. We eliminate this freedom by asking 
$t=0$ to be the first local maximum of $u_-$. Then $u_-'(0)=0$.
We note that $t=0$ is also a maximum of $v_0$, so $v'_0(0)=0$, and
 the unstable solution admits the approximation
\begin{equation}\label{Eq:outerapprox}
    u_-(t)=\varepsilon^2 v_0(\varepsilon t)+O(\varepsilon^4)
\end{equation}
which holds uniformly for $t \le t_0$ with any given $t_0$.
  In the leading term of this estimate,    $\varepsilon$ comes from the standard scaling.

The stable solution is  obtained using the reversible symmetry, $u_+(t)=u_-(-t)$. Since $v_0$ is even, $u_+$ admits approximation by the same function as $u_-$ but
  for $t\ge-t_0$.

 Then we consider the difference
 \[
 \delta(t)=u_+(t)-u_-(t).
 \]
 Since 	 $\delta(t)=u_-(-t)-u_-(t)$, it is odd and  we get  $\delta(0)=\delta''(0)=0$. The normalisation assumption
   $u_-'(0)=0$ implies that
		 $\delta'(0)=0$.
         Then,  in the case   $\delta'''(0)=0$, we can conclude that $u_+$ and $u_-$
         satisfy the same initial conditions at $t=0$ and,
 since the solution to Cauchy problem is unique,  
		 that $u_+(t)=u_-(t)$ for all $t$.
Hence we have a symmetric  soliton each time $\delta'''(0)$
vanishes.  This soliton is a single-bump one as it is approximated by  the unperturbed homoclinic solution $\varepsilon^2 v_0(\varepsilon t)$ for all real $t$.

A direct substitution shows that $\delta$
	satisfies the variational equation:
\begin{equation}\label{Eq:vareqn}
    	\delta^{(iv)}+(1-\varepsilon^2) \delta''-\varepsilon^2 \delta +2(1-\varepsilon^2)u\delta+2\gamma  (u\delta ''+u'\delta '+u''\delta)=0
\end{equation}
around $u=u_c=\frac12(u_++u_-)$.  The variational equation
is a linear homogeneous ODE of order four. The space of its
solution is spanned by four fundamental solutions. 
Let $\delta_1$ be a solution that satisfies  the Cauchy problem
for the variational equation:
$\delta_1(0)=\delta_1'(0)=\delta_1''(0)=0$ and $\delta_1'''(0)=1$. Then
there is a real constant $C=C(\varepsilon,\gamma)$ such that
\[
\delta(t)=C(\varepsilon,\gamma)\delta_1(t).
\]
We note that the existence of single-bump symmetric homoclinic solutions is equivalent to the vanishing of $C(\varepsilon,\gamma)$, i.e., it
 is a co-dimension one phenomenon.

\section{Inner equation and Stokes phenomenon\label{Se:inner}}

Determining the zeros of $C(\varepsilon,\gamma)$ is a non-trivial problem, since $C(\varepsilon,\gamma)$ is exponentially small in $\varepsilon$. The complex matching method  suggests studying the analytic continuation of the  solutions into the complex plane in order to obtain an  approximation for $C(\varepsilon,\gamma)$.

The leading order of this approximations comes from
the study of stable and unstable
solutions of the inner equation that is obtained from 
 \eqref{Eq:bifurcform} by setting $\varepsilon=0$,
\begin{equation}\label{Eq:inner}
u^{(iv)} + u'' + u^2 + \gamma\bigl({u'}^2 + 2u u''\bigr)=0.
\end{equation}
The equilibrium solution $u=0$ is a degenerate saddle-centre: the linearisation at the origin has two purely imaginary eigenvalues $\pm \mathrm i$ and a double zero eigenvalue. The inner equation 
admits two distinguished solutions, referred to as the stable and unstable solutions, which are uniquely characterised by their asymptotic behaviour
at the infinity:
\begin{equation}\label{Eq:upmasymp}
    u^\pm(x)=-\frac{6}{x^2}+O(x^{-4}), \qquad x\to\pm\infty,
\end{equation}
respectively, where we introduce yet a new time, $x = (\tau - i\pi)/\varepsilon$ (see Remark~\ref{rem:newTimeX} below).  Similarly to $u_\pm$, the stable and unstable solutions of \eqref{Eq:bifurcform}, the solutions $u^\pm$ belong to the stable and unstable manifolds of the equilibrium for equation~\eqref{Eq:inner},
 although the convergence rate is much slower.

\begin{rem}\label{rem:newTimeX}
The complex matching method
uses the observation that the approximation~\eqref{Eq:outerapprox} 
fails near the singularity of $v_0(\tau)$
located at $\tau=\mathrm i \pi$, and that in this region a better approximation is obtained using the inner equation. 
The method uses the change of time 
\begin{equation}
    \tau=\varepsilon t=\rmi \pi+\varepsilon x
\end{equation}
to study stable and unstable solutions in the region where  $|\tau-\mathrm i \pi |$ is small and $|x|$ is large.
The application of the standard scaling to 
the Laurent series $v_0(\tau)=-6(\tau-\rmi\pi)^{-2}+O(1)$
suggests to use \eqref{Eq:upmasymp} as an asymptotic
condition for the inner equation at infinity. 
Then it can be shown that $u^{\pm}(x)$ approximate  the stable and unstable solutions $u_\pm(t)$ 
with $t=\mathrm i\frac{\pi}{\varepsilon}+x$. 
 \end{rem}

Since equation~\eqref{Eq:inner} is reversible and
all its coefficients are real, the involution
\begin{equation}\label{Eq:involutionR}
(Ru)(x)=\overline{u(-\overline{x})}
\end{equation}
 maps a solution $u$ of~\eqref{Eq:inner} to another solution of the same equation. 
 Under this transformation, the domain of $u$ is   reflected in  the imaginary axis.
 Since the solutions $u^\pm$ are uniquely defined by the asymptotic condition \eqref{Eq:upmasymp} which is invariant
 under the involution $R$, we conclude that 
$
u^+ = Ru^-$,
or equivalently
\[
u^+(x)=\overline{u^-(-\overline{x})}.
\]
All points on the imaginary axis are fixed under the involution 
$x\mapsto -\overline x$.
 We consider the difference
 \[
\delta_u(x)=u^+(x)-u^-(x).
\]
on the imaginary axis for large negative values of~$\operatorname{Im}(x)$.

First we look for a formal solution of the inner equation 
\eqref{Eq:inner}.
The reversibility of the equation allows us
to use an ansatz that contains only even powers:
\begin{equation}\label{Eq:uf}
u_{\text{formal}}=\sum_{k\ge1}\frac{a_k}{x^{2k}}.
\end{equation}
After substituting the series into \eqref{Eq:inner} and collecting terms of equal order in $x$
we obtain a recurrent system for the coefficients $a_k$. The leading order gives $a_1=-6$
and 
equation 
\eqref{Eq:recursionak} provides an explicit recursive 
 formula for $a_k$ with $k>1$.

\begin{figure}
    \centering
\includegraphics[width=6cm]{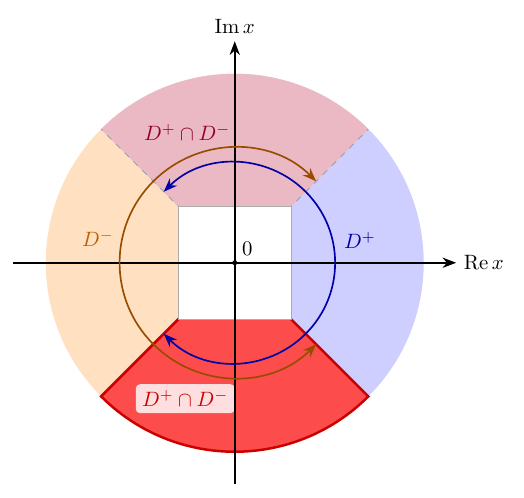}
    \caption{Sketch of the sectors $D^-$ and $D^+=R(D^-)$ where $u^\pm$ are defined.}
    \label{Fig:sectors}
\end{figure}

It can be shown that the formal solution provides an asymptotic expansion for both the stable and unstable solutions,
\begin{equation}\label{Eq:asymtcond}
u^\pm(x)\asymp u_{\text{formal}},
\end{equation}
as $x\to\infty$ in two large sectors (with opening angles greater than $\pi$) centred respectively around the positive and negative real semi-axes sketched in Fig.~\ref{Fig:sectors}.

We emphasize that, although the same formal series arises in both cases, the domains in which the partial sums of $u_{\text{formal}}$ approximate $u^\pm$ are different. On the intersection of these domains, the difference
$
\delta_u(x)=u^+(x)-u^-(x)
$
is beyond all algebraic orders, in the sense that $\delta_u(x)=O(x^{-n})$ for every $n\in\mathbb{N}$. This common domain of validity contains two sectorial neighbourhoods of infinity centred around the positive and negative imaginary semi-axes (see  Fig.~\ref{Fig:sectors}).

In order to obtain an estimate for $\delta_u$, we consider the variational equation
\begin{equation}\label{EQ:innervareqn}
\delta^{(iv)} + \delta'' + 2u\,\delta + 2\gamma\bigl(u\,\delta'' + u'\,\delta' + u''\,\delta\bigr)=0,
\end{equation}
which is obtained by linearising the inner equation~\eqref{Eq:inner} around a given solution $u$.
A direct substitution shows that $\delta_u$ satisfies~\eqref{EQ:innervareqn} with $u=u^c$, where
\[
u^c(x)=\tfrac12\bigl(u^+(x)+u^-(x)\bigr).
\]
Note that equation~\eqref{EQ:innervareqn} coincides with \eqref{Eq:vareqn}
with $\varepsilon=0$ and $u_c$ replaced by~$u^c$.

Equation~\eqref{EQ:innervareqn} is a fourth-order linear homogeneous ODE, so the space of its solutions is spanned by four linearly independent solutions. Taking into account the asymptotic behaviour of the coefficients in~\eqref{EQ:innervareqn}, one can show that these solutions may be chosen to behave, respectively, as $x^{-3}$, $x^4$, 
and $\mathrm e^{\pm \rmi x}$ as $x\to -\rmi\infty$ along $\rmi\mathbb{R}$. Of these four solutions, only two decay along the negative imaginary axis, and only one is $o(x^{-4})$.

To construct an approximation to the
fast decaying solution, we observe that the  variational equation \eqref{EQ:innervareqn} with $u=u_{\text{formal}}$ admits 
 a formal solution of the form
\begin{equation}\label{Eq:deltaformal}
\delta_{\text{formal}}=\mathrm e^{-\rmi  x}\sum_{k\ge 0} \frac{b_k}{x^{k}}.
\end{equation}
Substituting the formal series $u_{\text{formal}}$ and $\delta_{\text{formal}}$ into~\eqref{EQ:innervareqn} and collecting terms of equal order in $x$ yields a recursive system for the coefficients $b_k$. This system determines the coefficients uniquely up to a multiplicative constant, which we fix by setting $b_0=1$. Then the coefficients $b_k$ with even $k$ are real and the ones with odd $k$ are purely imaginary.

Since both $u^-$ and $u^+$ share the same asymptotic expansion on the intersection of their domains, it follows that $u^c \asymp u_{\text{formal}}$.
It can be shown that truncated formal series provide  approximation
to an analytic solutions of the variational equation. 
In particular,
 one can show that the variational equation with $u=u^c$ admits a unique solution,
\begin{equation}\label{Eq:deltaWKB}
\delta_{\text{WKB}}(x)\asymp \delta_{\text{formal}}
\end{equation}
as $x\to-\rmi \infty$ in a sector centred on the negative imaginary semi-axis.


%
%
We note that \eqref{Eq:asymtcond}
implies that $\delta_u$ decays
faster than any power along the imaginary axis. Therefore,
since $\delta_u$ is  a solution of~\eqref{EQ:innervareqn},
it must be proportional to $\delta_{\text{WKB}}$.

\begin{defn}[Stokes constant 1]
 The complex constant 
$S_0=S_0(\gamma)$
such that in a neighbourhood of infinity at the negative imaginary semi-axis
\begin{equation}\label{Eq:S0}
u^+(x)-u^-(x)=S_0 \, \delta_{\mathrm{WKB}}(x)
\end{equation}
is called a Stokes constant.
\end{defn}
Since 
$
\delta_{\mathrm{WKB}}(x)
= \mathrm e^
{-\mathrm i x} \left(1+O(x^{-1})\right)
$
we  get an equivalent description
for $S_0$.
\begin{prop}
The Stokes constant $S_0$ is equal to the following limit:
\begin{equation}\label{Eq:S0def}
S_0= \lim_{\Im x\to -\infty }\mathrm e^{\mathrm i x}\bigl(u^+(x)-u^-(x)\bigr).
\end{equation}
\end{prop}

\begin{prop}
    The Stokes constant $S_0$ is purely imaginary.
\end{prop}

This property follows from the reversibility and real-analytic symmetry of the inner equation. Indeed, 
 $\delta_u$ is antisymmetric with respect to involution $R$ defined in \eqref{Eq:involutionR},
\[
(R\delta_u)(x)=(Ru^+)(x)-(Ru^-)(x)=u^-(x)-u^+(x)=-\delta_u(x).
\]
On the other hand,
it can be checked using induction that $\overline{b}_k=(-1)^k b_k$ for all $k\ge0$. Therefore the formal expansion $\delta_{\text{formal}}$ is symmetric,
$R(\delta_{\text{formal}})=\delta_{\text{formal}}$. 
Since both
the equation and  asymptotic condition are invariant,
the uniquely defined solution $\delta_{\text{WKB}}$ is also symmetric,
%
i.e. \[\delta_{\mathrm{WKB}}(x)=R(\delta_{\mathrm{WKB}})(x).\]
Applying $R$ to the identity~\eqref{Eq:S0} therefore yields $S_0=-\overline{S_0}$, implying that the Stokes constant $S_0$ is purely imaginary.

\begin{rem}
If $S_0(\gamma)\ne0$ for some $\gamma$, 
 we get an example of the non-linear Stokes phenomenon. Informally we can say that 
the unstable solution $u^-$ changes its asymptotic behavior at infinity
from $u_{\text{formal}}$ to $u_{\text{formal}}-S_0(\gamma)\delta_{\text{formal}}$
when  $x$ crosses the Stokes line located on the negative imaginary semi-axis. 

We note that a special care is needed in the interpretation of this statement as it involves divergent series.  
Truncating the formal series at the  optimal order, i.e. at the
least term for $|x| = \sigma$, and denoting the corresponding partial
sum of the formal series by $u_n$, we visualize the Stokes phenomenon
in Fig.~\ref{fig:stokesphenomenon}: 
the difference $u^-(x) - u_n(x)$ transitions from monotone to oscillatory behaviour as $x$ crosses the imaginary axis 
following a line $\Im x=-\sigma$.
 
\begin{figure}[htp]
    \centering
    \includegraphics[width=0.48\linewidth]{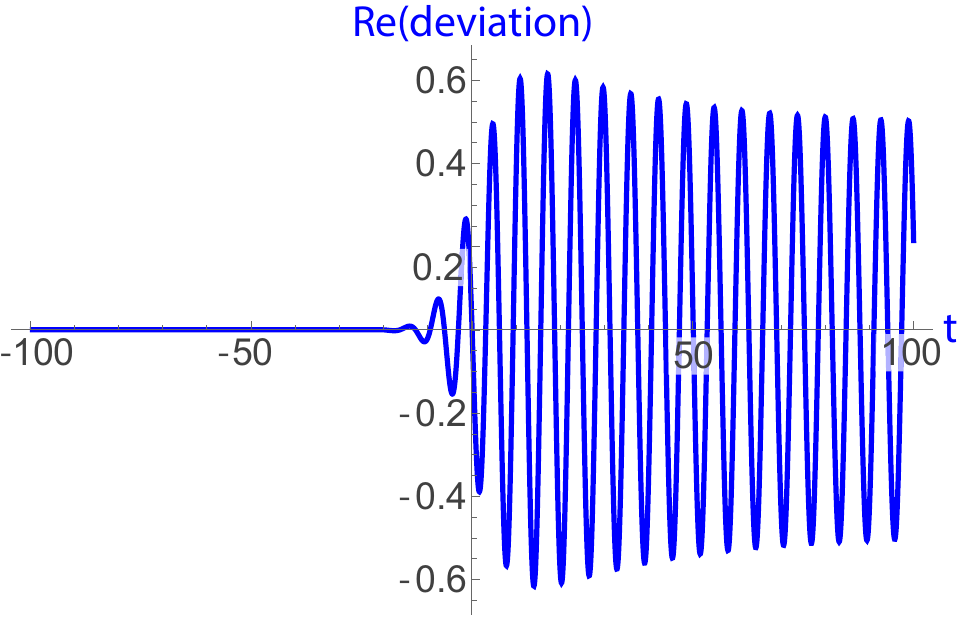}
    \includegraphics[width=0.48\linewidth]{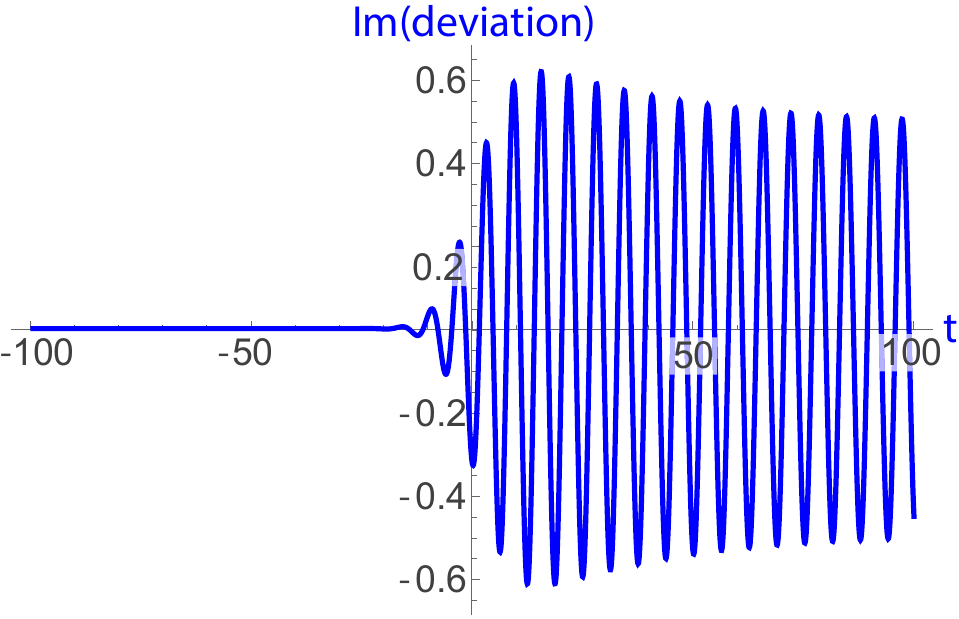}
    \caption{Stokes phenomenon: deviation of the unstable solution from the truncated asymptotic series. 
Plot of $\mathrm e^{\sigma}(u^-(x) - u_n(x))$ with $x = t - \mathrm i\sigma$ versus $t$, real part (left) and imaginary part (right). Parameters 
$\gamma=3$, $n=20$, $\sigma=40$.
    }
    \label{fig:stokesphenomenon}
\end{figure}
\end{rem}

\begin{prop}\label{prop:stokcons}
The Stokes constant $S_0=0$ if, and only if, the formal  series $u_{\text{formal}}$ converge.
\end{prop}

 Suppose that  $S_0(\gamma)=0$ for some $\gamma$. Then the equation \eqref{Eq:S0}
implies that  $u^+(x)-u^-(x)\equiv 0$
 near infinity on the negative imaginary semi-axis. The principle of analytic
 continuation implies that $u^+(x)=u^{-}(x)$ for all $x$ where the functions are defined. The union of the domains of $u^+$
 and $u^-$ includes  a complete neighbourhood of infinity
 and is not simply connected. 
The  real-analytic symmetry can be used to show
that $u^-=u^+$ is single-valued  in a neighbourhood of infinity. Thus its  Laurent series converge
and, due to the uniqueness,  coincide with   the formal solution $u_{\text{formal}}$.
Thus  the  series $u_{\text{formal}}$ have a finite radius of
convergence.

The reverse statement also holds: if 
the formal series converge, the sum satisfies the inner equation
and  both asymptotic conditions. Consequently the sum
 coincides both with $u^+$ and $u^-$. Thus  $u^+\equiv u^-$ and $S_0(\gamma)=0$.

\medskip
 
These arguments imply that the Stokes constant is related to
the growth rate of the coefficients in the  formal solution.

\begin{defn}[Stokes constant 2]
The Stokes constant $S_1=S_1(\gamma)$ is the limit
\begin{equation}
\label{Eq:S2}
S_1=\lim_{k\to\infty}\frac{(-1)^{k}a_k}{(2k-1)!}
\end{equation}
where $a_k$ are coefficients of the formal solution \eqref{Eq:uf}.
\end{defn}

\begin{rem}
Although it is not obvious, it can be shown that 
\[
S_0(\gamma)=-\pi{\mathrm  i} \,S_1(\gamma).
\] 
A sketch of the proof is in Appendix~\ref{Se:resurgent}.
This equality provides an efficient way for  computing the Stokes constant numerically. 
The equation~\eqref{Eq:S2} can also be  used to prove that the Stokes constant depends 
analytically on the parameter $\gamma$. Indeed, applying the induction to recursion \eqref{Eq:recursionak}, it is easy to check that 
 $a_k=a_k(\gamma)$ are polynomial of degree $k-1$
 with sign-alternating coefficients.  
In order to show that $S_1(\gamma)$ is analytic in $\gamma$ 
for $|\gamma|<\gamma_0$ it is sufficient to prove 
uniform convergence to the limit in \eqref{Eq:S2}.
\end{rem}

The Stokes constant is an analytic function of the real parameter $\gamma$ appearing in the inner equation~\eqref{Eq:inner}. Since $\Re S_0(\gamma)\equiv 0$, any change of sign of $\Im S_0(\gamma)$ implies the existence of a zero. By analyticity, all zeros are isolated and have finite multiplicity. Our numerical experiments suggest that all zeros are simple.

In the next section, we show that there are infinitely many zeros for $\gamma>1$ and derive an asymptotic description of their location in the limit $\gamma\to+\infty$.

\section{Asymptotics for the Stokes constant\label{Se:Sasympt}}

\subsection{Inner equation as a singularly perturbed problem}

In order to study the dependence of the Stokes constant on the parameter $\gamma$ 
for  $\gamma\gg1$ we introduce the new small parameter $\epsilon=\gamma^{-1/2}$.  The  substitution 
\begin{equation}\label{Eq:scaling1}
u(x)=\epsilon^{2}v(y)\quad\text{and}\quad y=\epsilon x
\end{equation}
transforms the inner equation~\eqref{Eq:inner}
into a singularly perturbed problem,
%
%
\begin{equation}
	\label{Eq:innersingular}
\epsilon^2 v^{(iv)} + v'' + v^2 +{v'}^2 +2vv''=0,
\end{equation}
where the small parameter appears in front of the highest  derivative. We stress that this equation depends on $\gamma=\epsilon^{-2}$
only and is completely independent
of  small parameter $\varepsilon$ introduced in \eqref{Eq:bifurcform}. From the viewpoint of
 singular perturbation theory, both parameters play a similar role due to the apparent similarity in the structure of equations \eqref{Eq:sing} and \eqref{Eq:innersingular}.

The unstable and stable solutions of the equation
\eqref{Eq:innersingular}  can be obtained by applying the scaling
 to
the corresponding solutions of equation~\eqref{Eq:inner}. 
In particular, after the substitution \eqref{Eq:scaling1},
the asymptotic condition~\eqref{Eq:upmasymp} takes the form
\begin{equation}\label{Eq:vpmasymp}
v^\pm(y)=-\frac{6}{y^2}+O(y^{-4})
\end{equation}
when $y\to\infty$ in sectors centred on the positive and negative
real semi-axes respectively.
We use  the method of complex
matching to approximate the difference
$$
\delta_v(y)=v^+(y)-v^-(y) 
$$
and to derive the asymptotic expression \eqref{Eq:stokesasympt} for the Stokes constant  in the limit of small $\epsilon$.

Equation \eqref{Eq:innersingular} admits a formal solution 
$v_{\mathrm{formal}}$, obtained from \eqref{Eq:uf} using the scaling 
\eqref{Eq:scaling1}. Applying  induction to the recurrence relations 
\eqref{Eq:recursionak}, one verifies that $a_k$ is a polynomial of degree 
$k-1$ in $\gamma=\epsilon^{-2}$. Hence $v_{\mathrm{formal}}$ can be 
written as a double series
\begin{equation}\label{Eq:vformal}
v_{\mathrm{formal}}=\sum_{\substack{k\ge0\\ 0\le j\le k-1}}
\frac{a_{kj}\,\epsilon^{2j}}{y^{2k}}
\end{equation}
where $a_{kj}\in\mathbb{R}$.
Since only even powers of $\epsilon$ occur, we look for a formal solution of 
\eqref{Eq:innersingular} in the form
\[
v=\sum_{k\ge0}\epsilon^{2k}v_k(y).
\]
Substitution into \eqref{Eq:innersingular} shows that the leading term 
$v_0$ satisfies the ODE
\begin{equation}\label{Eq:outer4inner}
v_0''+v_0^2+{v_0'}^2+2v_0v_0''=0,
\end{equation}
which is the equation obtained from \eqref{Eq:innersingular} in the limit 
$\epsilon\to0$. In view of the asymptotics \eqref{Eq:vpmasymp}, we select 
the solution satisfying
\begin{equation}\label{Eq:v0asymp1}
v_0(y)=-6y^{-2}+O(y^{-4}), \qquad y\to\infty.
\end{equation}

\subsection{Outer problem for the inner problem}
In the terminology of the complex matching, 
equation \eqref{Eq:outer4inner}
plays 
the role of an outer equation. It is a reversible non-linear
differential equation of the second order. The equation is integrable. Indeed, 
multiplying  \eqref{Eq:outer4inner} by $v_0'$,  integrating and setting the integration constant
to zero we get
\begin{equation}\label{Eq:v0}
(\tfrac12+v_0){v_0'}^2+\tfrac13 v_0^3=0.
\end{equation}
We look for a solution of this equation in the form
$v_0=-\tfrac12 q_0^2$. We fix the branch of $q_0$
by the asymptotic condition
\begin{equation}\label{Eq:q0ass}
q_0(y)=\frac{2\sqrt{3}}{y}+O(y^{-3}).
\end{equation}
After taking the square root of \eqref{Eq:v0} 
and choosing branches of the square root to be consistent with the asymptotic condition \eqref{Eq:q0ass} for large positive $y$, 
we get
\begin{equation}\label{Eq:q0diff}
\sqrt{1-q_0^2}\,q_0'=-\tfrac{1}{2\sqrt{3}}q_0^2.
\end{equation}
Dividing the equality by $q_0^2$ and integrating we get
an explicit expression for the inverse of $q_0$:
\begin{equation}\label{Eq:q0inverse}
\frac{\sqrt{1-q_0^2}}{q_0}+\arcsin q_0=\frac{y}{2\sqrt{3}},
\end{equation}
The inverse function theorem implies that this equation has a 
unique solution $q_0$ that maps a neighbouthood of infinity to a neighbourhood of zero. This solution
 is an odd analytic
function  and satisfies~\eqref{Eq:q0ass}.

In addition to the functions $u^\pm$, the definition of the Stokes constant \eqref{Eq:S0}
involves an oscilatory solution of the variational equation \eqref{EQ:innervareqn}. 
After applying the scaling \eqref{Eq:scaling1}, this equation takes the form
\begin{equation}
	\label{Eq:inner_variational}
	\epsilon^2 \eta^{(iv)} + \eta'' + 2v\eta 
	+ 2\bigl(v \eta'' + v'\eta' + v''\eta\bigr) = 0.
\end{equation}
In particular, 
the difference $\delta_v(y) = v^+(y) - v^-(y)$ satisfies 
\eqref{Eq:inner_variational} with $v = v_c(y) := \tfrac12\bigl(v^+(y) + v^-(y)\bigr)$. 
The WKB theory suggests to look for a solution 
in the form
\begin{equation}\label{Eq:yWKB}
\eta_{\text{WKB}}(y)=  \mathrm{e}^{- \mathrm i\,\epsilon^{-1}\phi_0(y)}\bigl(\eta_0(y)+\dots\bigr)
\end{equation}
where the dots stand for terms of higher orders in $\epsilon$.
Substituting this ansatz into the equation  \eqref{Eq:inner_variational},
taking into account that $v_c\sim  v_0$
and collecting the leading terms in $\epsilon$
we get (after simplifying the expressions) that 
\[
{\phi_0'}^2=1+2v_0,\qquad 
2(1+2v_0)\eta_0'+3v_0\eta_0=0.
\]
These equations can be solved explicitly in terms of $q_0$:
\begin{equation}\label{Eq:phi0eta0}
\phi_0=2\sqrt{3}
\left(q_0+q_0^{-1}\right),
\qquad
\eta_0=\left(1-q_0^2\right)^{-3/4}.
\end{equation}
Taking into account the properties of $q_0$
  we conclude that
$\phi_0$ and $\eta_0$ are both single-valued in a neighbourhood of infinity. Then \eqref{Eq:q0ass}
implies
\begin{equation}\label{Eq:phi0eta0Atinfty}
\phi_0(y)=y+O(y^{-1})\quad\text{and}\quad \eta_0(y)=1+O(y^{-2}).
\end{equation}
Finally, we note that both $\eta_{\text{WKB}}(y)$ and $ \delta_{WKB}(\epsilon^{-1}y)$ are solutions of \eqref{Eq:inner_variational} of the fastest decay along
the imaginary axis. 
Comparing \eqref{Eq:yWKB} with
 \eqref{Eq:deltaWKB}, 
we conclude that $\eta_{\text{WKB}}(y)= \delta_{WKB}(\epsilon^{-1}y)$ as their ratio is constant and converges to $1$
when  $y\to -\mathrm i \infty$.

\subsection{Singularities of the outer problem}

We note that the asymptotic  \eqref{Eq:v0asymp1} holds in a full 
neighbourhood of infinity; hence $v_0$ provides a common approximation 
to both $v^-$ and $v^+$. Complex matching suggests to consider  
analytic continuation of solutions into neighbourhoods of the singularities 
of $v_0$. Therefore we need information on the location and nature of the 
singularities of $v_0=-\frac12 q_0^2$. 

Equation \eqref{Eq:q0diff} is a non-linear first-order ODE  and the corresponding vector field grows as $O(q_0)$
at the infinity. Consequently,
 $q_0$ admits analytic continuation 
along any path in the complex $y$-plane such that the values of $q_0$ along this path avoid the singularities of 
the vector field located at $q_0=\pm1$. Using that
\[
\arcsin q = -\mathrm i \log\bigl(\mathrm i q + \sqrt{1 - q^2}\bigr),
\]
we conclude that equation \eqref{Eq:q0inverse} implies that the singular values of the  function $y(q_0)$ are located at the points
\[
y_k = \pi\sqrt{3}\,(2k+1), \qquad k\in\mathbb Z.
\]
Consequently, $q_0(y)$ admits analytic continuation 
along any path which avoids $y_k$.

For the principal branch of the inverse sine function, the left-hand side of 
\eqref{Eq:q0inverse} is a monotone function mapping $(0,1)$ onto 
$(y_*,+\infty)$, where $y_* := y_0 = \pi\sqrt{3}$. Consequently, 
$q_0$ is analytic on $(y_*,+\infty)$. Since $q_0$ is odd, it is also 
analytic on $(-\infty,-y_*)$. Hence the principal branch of $q_0$ 
extends analytically to $\mathbb{C}\setminus[-y_*,y_*]$, i.e.\ the 
complex plane cut along the interval $[-y_*,y_*]$.

In order to study $q_0(y)$ near the singularity at $y=y_*$, we recall the
 trigonometric identity that holds for all $q\in(0,1)$:
\[
\arcsin q = \arctan\frac{q}{\sqrt{1-q^2}}=\frac\pi2-\arctan\frac{\sqrt{1-q^2}}{q}.
\]
Then 
\eqref{Eq:q0inverse} implies that for $q_0\in(0,1)$
\[
\frac{\sqrt{1-q_0^2}}{q_0}-
\arctan\frac{\sqrt{1-q_0^2}}{q_0}
=\frac{y}{2\sqrt{3}}-\frac\pi2.
\]
Recalling the Maclaurin series 
\(
\arctan x = x - \frac{x^3}{3}+\frac{x^5}{5}+\dots
\)
and applying the inverse function theorem we get that
\[
\frac{\sqrt{1-q_0^2}}{q_0}= 3^{1/3}\left(\frac{y}{2\sqrt{3}}-\frac\pi2
\right)^{1/3}+O(y-y_*)
\]
where the remainder is a converging  series in odd powers of  $\left(y-y_*
\right)^{1/3}$. Taking the square  and resolving the result with respect to $q_0$ we get that
\[
q_0(y)=1-\tfrac12 \left(\tfrac{3}{4}\right)^{1/3}\left(y -y_*\right)^{2/3}+O(y-y_*)^{4/3} 
\]
where the remainder is a converging power series in $\left(y -y_*\right)^{2/3}$. Therefore $y=y_*$ is a rational branching point
with $q_0(y_*)=1$. 

Finally we note that the uniqueness of the initial value problem for the differential equation \eqref{Eq:q0diff} implies that $q_0$ does not vanish. 
We conclude that the functions $v_0$, $\phi_0$ and $\eta_0$
are  analytic on $\mathbb C\setminus[-y_*,y_*]$. Moreover
\begin{align}
\label{Eq:v0sing}
v_0(y)&=-\tfrac12q_0(y)^2=-\tfrac12+\tfrac12\left(\tfrac{3}{4}\right)^{1/3}(y-y_*)^{2/3}+O(y-y_*)^{4/3}   \\
\phi_0(y)&=2\sqrt{3}
\left(q_0(y)+q_0(y)^{-1}\right)
\nonumber
\\
&
=4\sqrt3+\left(\tfrac{3}{4}\right)^{7/6}(y-y_*)^{4/3}+O(y-y_*)^2
,
\label{Eq:phi0sing}
\\
\label{Eq:eta0sing}
\eta_0(y)&=\left(1-q_0^2(y)\right)^{-3/4}=\left(\tfrac{4}{3}\right)^{1/4}(y-y_*)^{-1/2}\left(1+O(y-y_*)^{2/3}\right)
.
\end{align}
The value of $\phi_0(y_*)=4\sqrt{3}$ will play an important role
in the subsequent analysis.

\subsection{Inner problem for the inner problem\label{Se:inner4inner}}

We have established that $y_*$ is a branching point of $v_0$ and
$v_0(y_*)=-\frac12$. In order to study $v^+ $ in a neighbourhood of this point, 
we make the substitution
\begin{equation}\label{Eq:subst2}
v(y)=-\tfrac12+\epsilon^{1/2} w(z),\qquad  y=y_*+\epsilon^{3/4} z.
\end{equation}
In the new variables the equation~\eqref{Eq:innersingular}
takes the form
\begin{equation}
  w^{(iv)} 
+ \epsilon^{1/2} \left(-\tfrac12+\epsilon^{1/2} w\right)^2 +{ w'}^2 
+2 w  w''
=0.
\end{equation}
We note that this equation represents a regular perturbation problem
as the derivative of order four is no longer multiplied by the small parameter. 
For $\epsilon=0$ we get the inner equation for the inner equation:
\begin{equation}\label{Eq:inner4inner}
w_0^{(iv)}  +{w_0'}^2  + 2 w_0  w_0'' =0.
\end{equation}
We look for solutions $w_0^\pm(z)$ of \eqref{Eq:inner4inner} whose behaviour 
at large  $z$ matches that of $v_0(y)$ near $y_*$, as described by 
\eqref{Eq:v0sing}. This asymptotic form suggests a formal solution of the form
\[
w_{\mathrm{formal}} = z^{2/3}\sum_{k\ge0}\frac{c_k}{z^{8k/3}}.
\]
Substitution of this series into \eqref{Eq:inner4inner} shows that the 
coefficients $c_k$ for $k\ge1$ are uniquely determined once $c_0$ is fixed. 
Comparison with \eqref{Eq:v0sing} then leads to the choice
$
c_0 = \tfrac12\left(\tfrac{3}{4}\right)^{1/3}
$.

The variational equation associated with \eqref{Eq:inner4inner} is given by
\begin{equation}\label{Eq:inner4innerVareqn}
\xi^{(iv)} + 2\bigl(w_0\xi'' + w_0'\xi' + w_0''\xi\bigr) = 0.
\end{equation}
This is a linear homogeneous ODE of order four and, for $w_0 = w_{\mathrm{formal}}$,  the space of 
its solutions is spanned by four independent solutions that 
respectively behave as $z^{-1/3}$, $z^{2/3}$ and $z^{-1/2}\mathrm e^{\pm \rmi c z^{4/3}}$ for large $z$.
The formal exponentially decaying solution
has the form
\begin{equation}\label{Eq:xiformal}
\xi_{\mathrm{formal}} = \exp\!\left(-\mathrm i\left(\tfrac{3}{4}\right)^{7/6}z^{4/3}\right)
z^{-1/2}\sum_{k\ge0}\frac{d_k}{z^{4k/3}}.
\end{equation}
The coefficients $d_k$ for $k\ge1$ are uniquely determined by the equation 
once $d_0$ is fixed. We set $d_0 = 1$. The exponential factor decays in the 
sector $\arg z \in (\alpha_-,\alpha_+)$, where $\alpha_+ = 0$ and 
$\alpha_- = -\tfrac{3\pi}{4}$.

The equation \eqref{Eq:inner4inner} admits two analytic solutions 
$w_0^\pm$ such that
\begin{equation}\label{Eq:w0asymp}
w_0^\pm(z)\asymp w_{\mathrm{formal}}(z),
\end{equation}
as $z\to \mathrm e^{\mathrm i\alpha_\pm}\infty$, respectively. 
Inside the sector between these directions, the difference
\[
\delta_w(z) = w_0^+(z) - w_0^-(z)
\]
is superalgebraically small, i.e.\ it decays faster than $z^{-n}$ for any 
$n>0$. A direct computation shows that $\delta_w$ satisfies the variational 
equation \eqref{Eq:inner4innerVareqn} with $w_0$ replaced by 
$w_c := \tfrac12(w_0^+ + w_0^-)$. 

Since \eqref{Eq:inner4innerVareqn} with $w_0=w_c$ is a linear ODE of order four, its solution 
space is four-dimensional. Among its solutions, only one,
\[
\xi_{\mathrm{WKB}}(z) \asymp \xi_{\mathrm{formal}}(z),
\]
decays faster than $z^{-1}$ along the bisecting direction of the sector. 
It follows that there exists a complex constant $K_0$ such that
\begin{equation}\label{Eq:w0Stokes}
w_0^+(z) - w_0^-(z) = \delta_w(z) = K_0\,\xi_{\mathrm{WKB}}(z).
\end{equation}
The equation~\eqref{Eq:inner4innerVareqn} has a 
symmetry which involves the reflection of the $z$-plane in the line $\arg z=\frac{3\pi}{8}$, the bi-section of the sector, and the reflection of the $w_0$-plane in
the line $\arg w_0=\frac{3\pi}4$: if $w_0(z)$ is a solution then $\mathrm e^{-\rmi 3\pi/2 }\overline{w_0(\mathrm e^{-\rmi 3\pi/4 }\overline z)}$ also is.
This symmetry swaps $w_0^+$ and $w_0^-$ and multiplies  $\xi_{\text{formal}}$ and, consequently, $\xi_{\text{WKB}}$ by $\mathrm e^{\rmi \frac{15\pi}{8}}$.
Applying this symmetry to \eqref{Eq:w0Stokes} 
we conclude that $K_0\mathrm e^{-\mathrm i\frac\pi{16}}$ is real.
Evaluating $\delta_w$ and $\xi_{\text{WKB}}$ numerically we find that
\[
|K_0|=0.251346943975660732\dots
\]

\subsection{Complex matching for the inner problem}

Now we are ready to apply the complex matching 
to estimate  the Stokes constant 
for equation \eqref{Eq:innersingular} with small $\epsilon>0$.
First we note that the functions $v^+$ and $v_0$ have the same 
asymptotic behaviour 
(given by \eqref{Eq:vpmasymp} and
\eqref{Eq:v0asymp1} respectively) 
at infinity in a large sector  centered
around the positive semi-axis. 
Consequently, the difference
\[
\delta^+(y)=v^+(y)-v_0(y)
\]
vanishes as $O(y^{-4})$ in the limit $\Re(y)\to+\infty$.

We recall that $v^+$ is an exact solution of  \eqref{Eq:innersingular} 
and while $v_0$ is an approximate one.
Since the growth of small errors is
governed by the variational equation
it is natural to study $\delta^+$
along a line 
\[
\ell_\sigma=\bigl\{y\in\C:\Im \phi_0(y)=\sigma\bigr\}
\]
as the solutions of the variational equation are
oscillatory and do not grow exponentially along this line.
Since $\phi_0(y)\approx y$ for large $y$, it is easy to see that
  $\ell_\sigma$ with a large $\sigma$
is approximately horizontal.
Solving  equation \eqref{Eq:phi0eta0} with respect to $\phi_0$
we get
\[
 q_0 = \frac{\phi_0}{\sqrt{48}}\left(1 - \sqrt{1-\frac{48}{\phi_0^2}}\right).
\]
Substituting this expression  into \eqref{Eq:q0inverse}
and simplifying the result we get
the parametrization for $\ell_\sigma$, 
\[
y= \phi_0\sqrt{\tfrac12+\tfrac12\sqrt{1-\tfrac{48}{\phi_0
   ^2}}}
   \left(1-\tfrac{48}{\phi_0 ^2}\right)^{1/4}
   +2 \sqrt{3} \arcsin
   \left(
   \tfrac{\phi_0}{\sqrt{48}}\left(1 - \sqrt{1-\tfrac{48}{\phi_0^2}}\right)\right)
.
\]
Let $\sigma = -\sigma_0 \epsilon$ for some $\sigma_0 > 0$. The corresponding 
curve lies in the lower half-plane. In the limit $\epsilon \to 0$, it 
converges to a contour consisting of two rays, $(-\infty,-\pi\sqrt{3}]$ and 
$[\pi\sqrt{3},+\infty)$, joined by an arc connecting $-\pi\sqrt{3}$ and 
$\pi\sqrt{3}$.

%
%
%
%
%
%
Let $y_k \in \ell_\sigma$, $k=1,2$, be the intersection points of 
$\ell_\sigma$ with the circle $|y-y_*|=\epsilon^{1/2}$. These points lie in 
a neighbourhood of $y_*$ on opposite sides of the Stokes line defined by $\Re\phi_0(y)=\Re\phi_0(y_*)$, so that
\[
\Re \phi_0(y_1) < \Re \phi_0(y_*) < \Re \phi_0(y_2).
\]
We use these points for complex matching. By \eqref{Eq:subst2}, the circle 
$|y-y_*|=\epsilon^{1/2}$ corresponds to $|z|=\epsilon^{-1/4}$. Hence, it 
represents a neighbourhood of $y_*$ in the $y$-plane and a neighbourhood of 
infinity in the $z$-plane. The scaling is chosen so that the errors of the 
inner and outer approximations are of the same order.

For $y \in \ell_\sigma$ to the right of $y_2$, the function $v^+(y)$ remains 
sufficiently far from the singular value $-\tfrac12$, so that the outer 
approximation is valid and $v^+(y) \approx v_0(y)$. In a neighbourhood of 
$y_2$, we  compare \eqref{Eq:v0sing} with the leading term in 
\eqref{Eq:w0asymp} for $w_0^+$ and conclude that
\[
v_0(y)\approx -\tfrac12+\tfrac12\left(\tfrac{3}{4}\right)^{1/3} (y-y_*)^{2/3} =-\tfrac12+\epsilon^{1/2}\tfrac12\left(\tfrac{3}{4}\right)^{1/3} z^{2/3}
\approx
-\tfrac12+\epsilon^{1/2} w_0^+(z).
\]
Consequently near $y=y_2$
\[
v^+(y)\approx v_0(y)\approx -\tfrac12+ \epsilon^{1/2} w_0^+(z).
\]
When $y$ lies on $\ell_\sigma$ between $y_1$ and $y_2$, the solution 
$v^+(y)$ remains close to $-\tfrac12$, and the inner approximation provides 
a more accurate description than the outer one. In a neighbourhood of $y_2$, 
we  use \eqref{Eq:w0Stokes} to switch from the stable solution $w_0^+$ 
to the unstable solution $w_0^-$:
\[
v^+(y)
\approx 
-\tfrac12+ \epsilon^{1/2}
\bigl(w_0^-(z)+K_0\xi_{\text{WKB}}(z)\bigr).
\]
Then 
 we compare \eqref{Eq:v0sing} and the first term of
\eqref{Eq:w0asymp} for $w_0^-$ and conclude that 
\[
v_0(y)\approx -\tfrac12+\tfrac12\left(\tfrac{3}{4}\right)^{1/3} (y-y_*)^{2/3} =-\tfrac12+\epsilon^{1/2}\tfrac12\left(\tfrac{3}{4}\right)^{1/3} z^{2/3}
\approx
-\tfrac12+\epsilon^{1/2} w_0^-(z).
\]
Then we compare 
 the approximation of $\xi_{\text{WKB}}$ 
 by the first term of \eqref{Eq:xiformal}
 and the approximation of 
$\eta_{\text{WKB}}$ by \eqref{Eq:yWKB},  \eqref{Eq:v0sing}
and \eqref{Eq:eta0sing}. We get that
\[
\begin{split}
\xi_{\text{WKB}}(z)&\approx
\exp\left(-\mathrm i\left(\tfrac{3}{4}\right)^{7/6}z^{4/3}\right)
z^{-1/2}
\\&=
\exp\left(-\mathrm i\left(\tfrac{3}{4}\right)^{7/6}\epsilon^{-1}(y-y_*)^{4/3}\right)\epsilon^{3/8}(y-y_*)^{-1/2}
\\
&\approx
\left(\tfrac{3}{4}\right)^{1/4}\epsilon^{3/8}
\exp\left(-\mathrm i \frac{-4\sqrt3+\phi_0(y)}{\epsilon}
\right)\eta_0(y).
\end{split}
\]
Consequently,  near $y_2$
\[
v^+(y)
\approx 
v_0(y)+
 \epsilon^{7/8}K_0
\left(\tfrac{3}{4}\right)^{1/4}
\exp\left(\mathrm i \frac{4\sqrt3}{\epsilon}\right)
\eta _{\text{WKB}}(y)
\]
Taking into account that small errors in initial conditions
follow solutions of the variational equation we get that 
this approximations extends to a neighbourhood of 
intersection of $\ell_\sigma$ with the imaginary axis.
Then we  apply the reversor to get the unstable
solution. For a point $y$ on the imaginary axis $\overline y=-y$ and,
consequently,
\[
\begin{split}
v^-(y)&=\overline{v^+(-\overline{y})}=
\overline{v^+(y)}
\\&
\approx 
v_0(y)+
 \epsilon^{7/8}\overline{K}_0
\left(\tfrac{3}{4}\right)^{1/4}
\exp\left(-\mathrm i \frac{4\sqrt3}{\epsilon}\right)
\eta _{\text{WKB}}(y).
\end{split}
\]
Taking the difference of the last two approximations and using 
$K_0 = \mathrm e^{\mathrm i\frac{\pi}{16}} |K_0|$, we obtain that on the 
imaginary axis
\[
v^+(y)-v^-(y)
\approx
\epsilon^{7/8}|K_0|
\left(\tfrac{3}{4}\right)^{1/4}
2\mathrm i \sin \left(\frac{\pi}{16}+ \frac{4\sqrt3}{\epsilon}\right)
\eta _{\text{WKB}}(y).
\]
Comparing the last estimate with the definition \eqref{Eq:S0} of the Stokes 
constant, and recalling that $u(x)=\epsilon^2 v(y)$ and 
$\xi_{\mathrm{WKB}}(x)=\eta_{\mathrm{WKB}}(y)$, we obtain
\[
S_0 \approx 
\epsilon^{2+7/8}|K_0|
\left(\tfrac{3}{4}\right)^{1/4}
2\mathrm i \sin \left(\frac{\pi}{16}+ \frac{4\sqrt3}{\epsilon}\right)
%
\]
Finally, recalling that $\epsilon=\gamma^{-1/2}$, we 
recover the asymptotic expression \eqref{Eq:stokesasympt}.

\section{Numerical methods for evaluation of the Stokes constant\label{Se:numerics}}

\subsection{Methods based on  analysis of the formal solution}

The most efficient method for computing the Stokes constant is based on the analysis of  the coefficients $a_k$ of the formal solution~\eqref{Eq:uf} to the inner equation~\eqref{Eq:inner}. 
These coefficients are defined by the recurrence relation with initial condition $a_1 = -6$ and
\begin{equation}\label{Eq:recursionak}
\begin{split}
a_k &= \frac{-1}{
  2 (2 k^2+k-6 )}
   \biggl(
(2 k - 2) (2 k - 1) (2 k) (2 k + 1) a_{k-1}
\\&\qquad\qquad +\sum_{j=2}^{k-1} a_{j}a_{k+1-j}
 +4\gamma 
 \sum_{j=1}^{k-1} 
 j(k+1+j)a_{j}a_{k-j}
 \biggr)
 \quad \text{for $k\ge 2$},
 \end{split}
\end{equation}
where the first of the sums at the rhs of~\eqref{Eq:recursionak} must be omitted for $k = 2$.
The induction in $k$ implies that this equation defines $a_k$
which are polynomial of degree $k-1$ in $\gamma$ with sign-alternating coefficients, and satisfy $a_k(5/4)=0$ for all $k>1$. 

\begin{figure}[htb]
    \centering
    \includegraphics[width=0.48\linewidth]{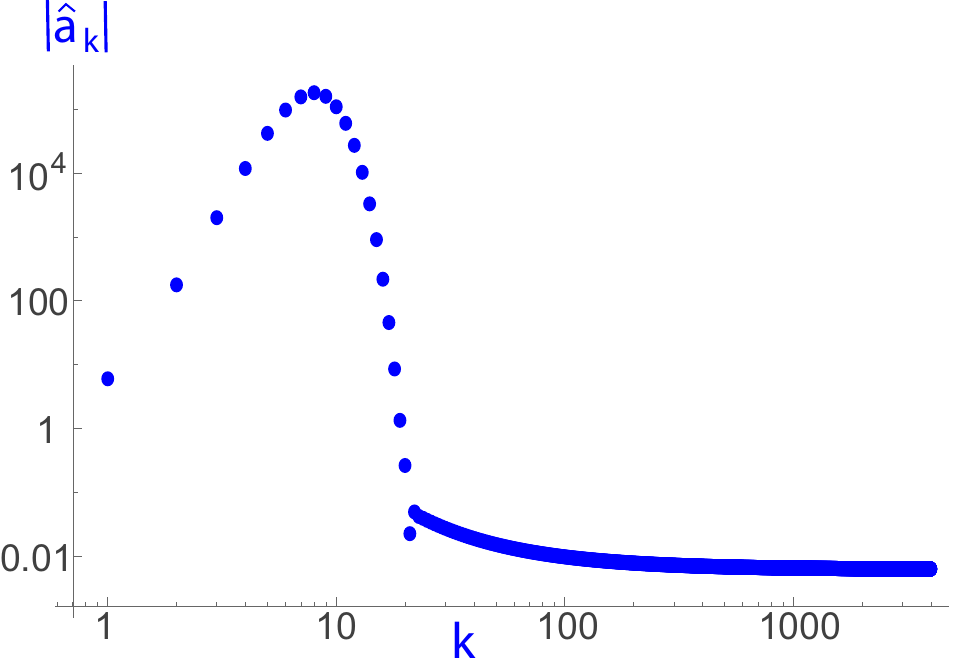}
    \includegraphics[width=0.48\linewidth]{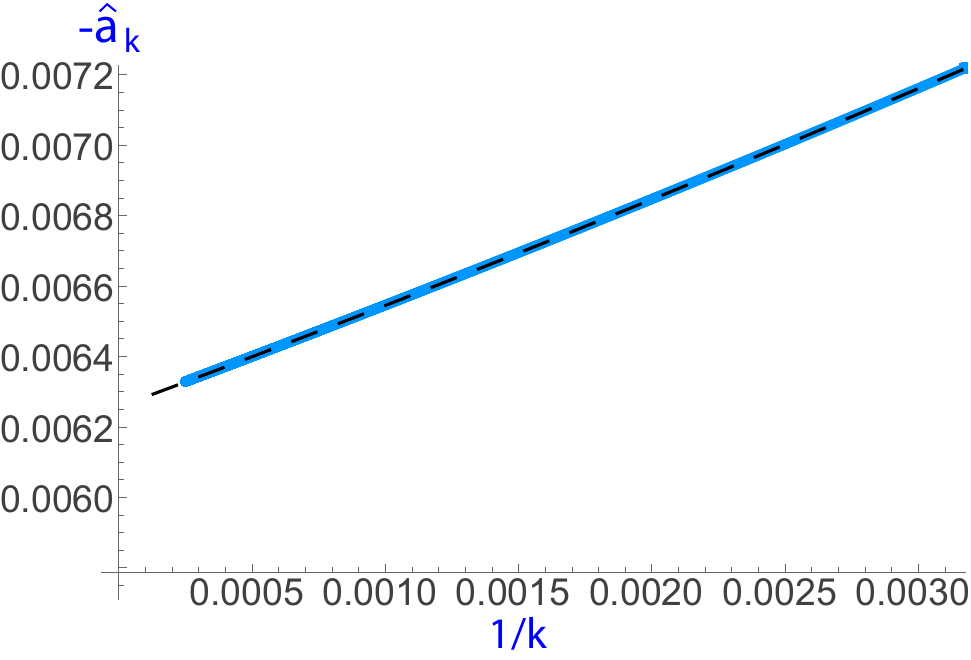}
    \caption{
    Normalized coefficients of the formal solution for $\gamma = 4^2$: 
$|\hat a_k|$ for $1 \le k \le 4000$ in log--log scale (left) and 
$-\hat a_k$ versus $k^{-1}$ (right); the thin discontinuous line shows an extrapolation to zero.}
    \label{fig:ak}
\end{figure}

For large $k$, the coefficients exhibit Gevrey-one growth with a sign-alternating tail.
We introduce the normalized coefficients
\begin{equation}
    \hat a_k = \frac{(-1)^k a_k}{(2k-1)!},
\end{equation}
and aim to compute the Stokes constant as the limit
\[
S_1 = \lim_{k \to \infty} \hat a_k.
\]
A typical behaviour of $\hat a_k$ for $\gamma > 1$ is shown in Fig.~\ref{fig:ak} (left), where we take $\gamma = 4^2$ as an example.
The tail of the sequence is monotone and decreases in absolute value. 
Plotting $\hat a_k$ as a function of $k^{-1}$ (Fig.~\ref{fig:ak}, right), we observe apparent convergence to a finite limit. 
The last of the computed terms of the series can be taken as an approximate value of the limit.
Much higher accuracy can be achieved  by extrapolation in powers of $k^{-1}$.

To this end, we proceed as follows. We identify two consecutive intervals $[k_0, k_0 + m]$ and $[k_0 + m, k_0 + 2m]$ on which the sequence $\hat a_k$ is monotone. 
On each interval, we construct interpolating polynomials $p_{m0}$ and $p_{m1}$ of degree $m$ using the data $(k^{-1}, \hat a_k)$. 
We  take $p_{m1}(0)$ as an estimate of $S_1$, provided that
the predicted values are in good agreement:
\begin{equation}\label{Eq:extrapolerror}
|p_{m1}(0) - p_{m0}(0)| < 10^{-50}.
\end{equation}
If this criterion is not satisfied, we increase $k_0$ and repeat the procedure.
We use the difference 
\eqref{Eq:extrapolerror} as a proxy for the error in evaluation of $S_1$.

The computation time depends on three factors:
the starting point $k_0$ that determines the number of coefficients $\hat a_k$ to be computed,
the number of correct digits  $N$
in numerical values of $\hat a_k$
and 
the degree $m$ of the interpolating polynomial.

As an example, we consider $\gamma = 4^2$. Using the tail of $4000$ coefficients and an interpolating polynomial of degree $80$, we obtain an estimate of $S_1(4^2)$ with an error below $10^{-180}$. 
This value is then used as a reference in subsequent, less accurate computations.

\begin{figure}[htb]
    \centering
    \includegraphics[width=0.6\linewidth]{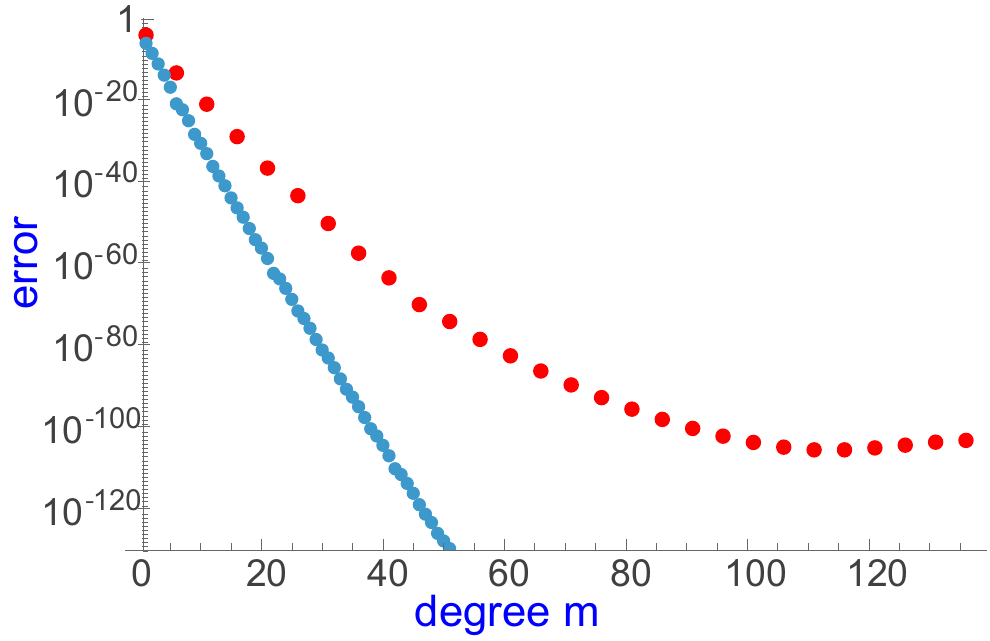}
    \caption{Accuracy of the extrapolation to zero as a function of the extrapolating polynomial degree $m$ using the tail of the first $4\,000$
    coefficients (blue) and $400$ coefficients (red).
    }
    \label{fig:extrapolation}
\end{figure}

We perform two series of experiments based on the tails of the first $400$ and $4000$ coefficients $\hat a_k$, respectively, and compare the  estimates for $S_1(\gamma)$ obtained using interpolating polynomials of varying degree. 
The deviation from the reference value for $S_1(4^2)$ is shown in Fig.~\ref{fig:extrapolation}. 
We observe that the error initially decreases rapidly as the polynomial degree $m$ increases, but eventually begins to grow.

 In particular, using the tail of the first $400$ coefficients
and an interpolating polynomial of degree $m = 100$, we obtain an estimate of the Stokes constant which deviates by less than  $10^{-100}$ from the reference value. 
The experiment confirms the numerical stability of the extrapolation procedure, with estimates of $S_1$ based on different intervals and different  degrees of extrapolating polynomials agreeing up to over 100 decimal digits.

In this experiment we also check 
  that~\eqref{Eq:extrapolerror} provides a reliable indicator of the extrapolation error.

We emphasize that the use of high-precision arithmetic is essential for the reliable implementation of this method, for two closely related reasons. 
First, repeated applications of the recurrence relation~\eqref{Eq:recursionak} lead to a significant loss of numerical precision. 
For instance, comparing computations performed with working precision $N = 500$ and $N = 1000$, we observe that approximately 180 digits are lost in the evaluation of $\hat a_k$ for $k \le  1000$ and $\gamma = 100^2$. 
Second, the extrapolation  using a high-degree polynomial is numerically stable only when the input data is computed with sufficiently high accuracy.

In the majority of our computation we use $m=50$
to find $k_0$ suitable for  
numerically reliable extrapolation.
Then we use $m=100$ to get a high precision
approximation of $S_1$.

\begin{figure}[htb]
     \centering
    \includegraphics[width=0.48\linewidth]{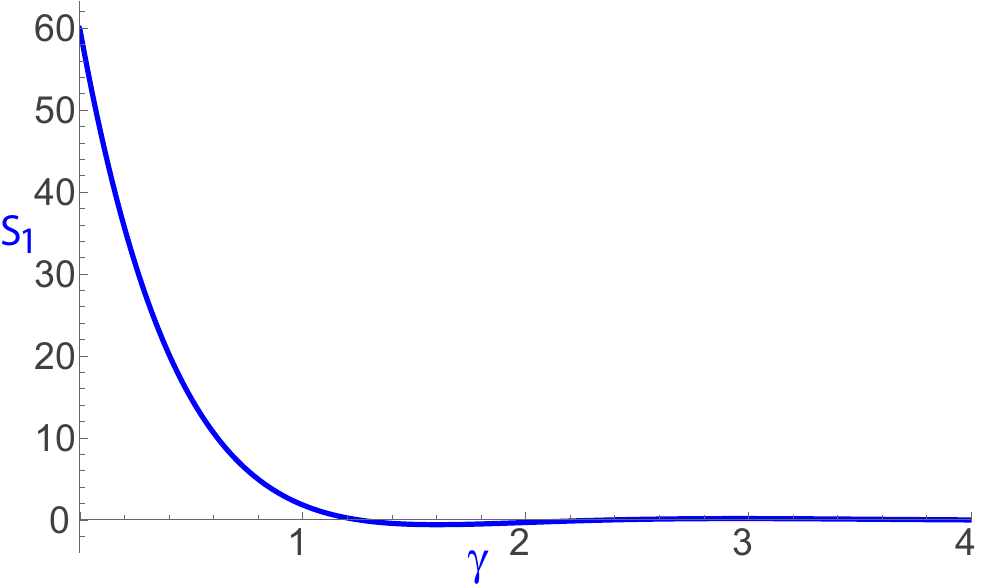}
\includegraphics[width=0.48\linewidth]{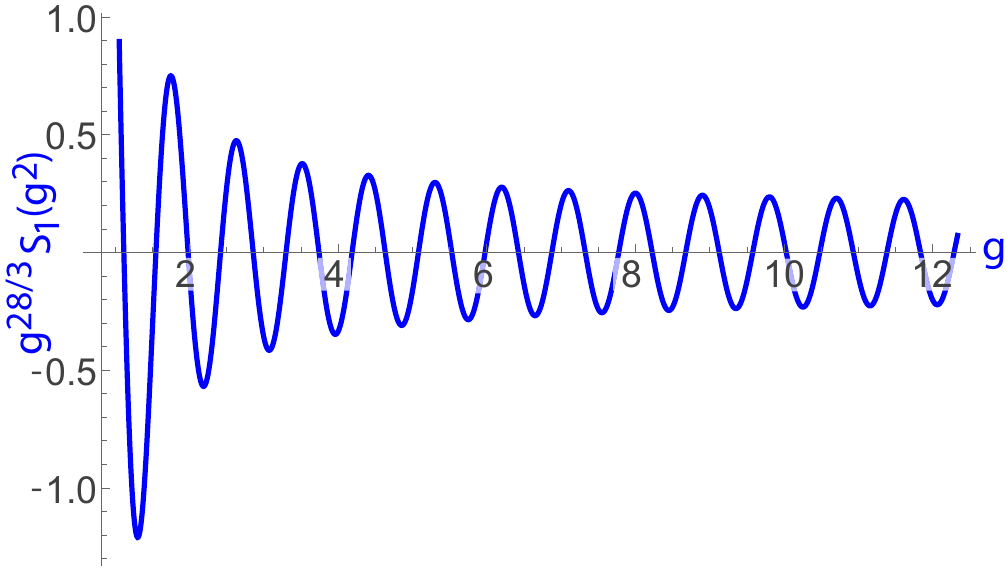}
     \caption{Stokes constant $S_1$ as a function of $\gamma$ (left), and the normalized value $g^{23/8} S_1(g^2)$ as a function of $g=\gamma^{1/2}$ on the interval which contains the first 26 zeroes of $S_1$ (right).
     }
     \label{fig:S1gamma}
 \end{figure}

The results of these computations are presented in Fig.~\ref{fig:S1gamma}. They show that $S_1(\gamma)$ is positive and monotonically decreasing for $\gamma < \frac{5}{4}$. It has its first zero at $\gamma=5/4$.
 Beyond this point, the behaviour becomes oscillatory.

To illustrate the oscillatory  behavior, we adopt the scaling suggested by the asymptotic estimate~\eqref{Eq:stokesasympt} and plot $g^{23/8} S_1(g^2)$ as a function of $g = \gamma^{1/2}$. For large $g$, the data clearly indicate convergence to a periodic function. A more accurate analysis of maximal and minimal
values confirms correctness of our asymptotic formula for
the Stokes constant.

\subsection{Zeroes of the Stokes constant}

\begin{figure}[htb]
    \centering
    \includegraphics[width=0.6\linewidth]{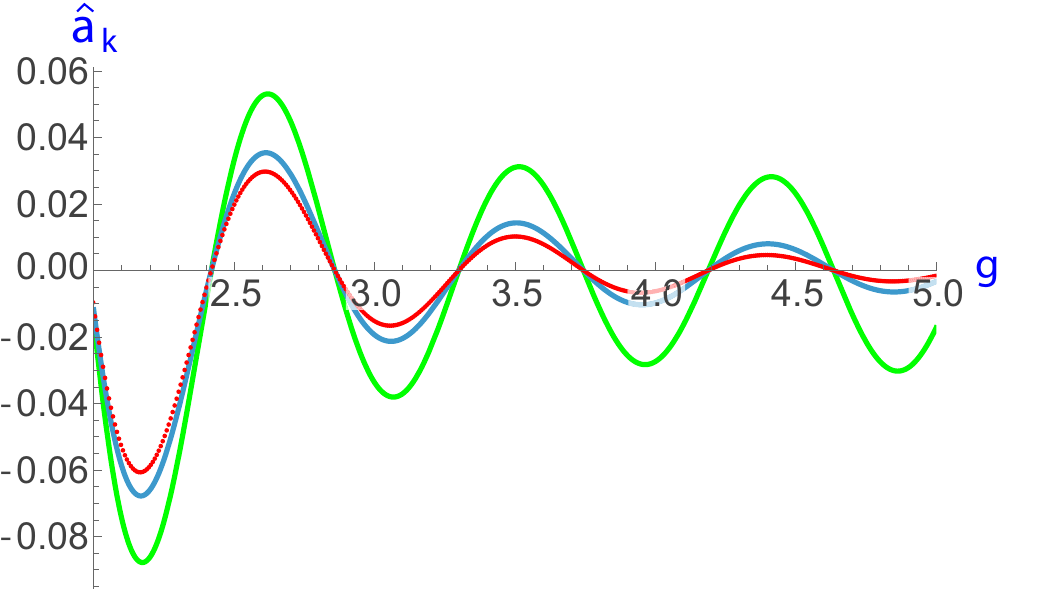}
    \caption{Normalized coefficients $ \hat{a}_k$ as functions of $g=\gamma^{1/2}$ for $k = 30$ (green) and $k = 100$ (blue), and the limiting curve $S_1(g^2)$ (red). }
    \label{fig:ak-g}
\end{figure}

We use the behavior of the normalized coefficients $\hat a_k$ to study the dependence 
of the Stokes constant $S_1$ on $\gamma$.
It is easy to see that $\hat a_k$ is 
polynomial of degree $k-1$. For example,
\[
\hat a_1=6,\quad
\hat a_2=-3 (4 \gamma -5),\quad 
\hat a_3=
\frac{3}{20} \left(64 \gamma ^2-204 \gamma +155\right).
\]
The polynomials with $k>1$ share the common root at $\gamma=5/4$.
For this special value of $\gamma$, $u(x)=-6x^{-2}$ solves the inner equation~\eqref{Eq:inner} exactly.
Figure~\ref{fig:ak-g}
illustrates convergence of $a_k$ to the limit.

We observe that the zeros of the Stokes constant $S_1(\gamma)$ can be located without relying on the extrapolation procedure. 
Indeed, suppose that $S_1(\gamma) \neq 0$ for some $\gamma$. Then there exists $k_0(\gamma)$ such that
\[
\operatorname{sign} S_1(\gamma) = \operatorname{sign} \hat a_k(\gamma)
\quad \text{for all } k \ge k_0(\gamma).
\]
Therefore, if
\[
\hat a_k(\gamma_1)\,\hat a_k(\gamma_2) < 0
\]
for some $k > \min\bigl(k_0(\gamma_1), k_0(\gamma_2)\bigr)$, 
then the interval $(\gamma_1, \gamma_2)$ contains a zero of $S_1(\gamma)$. 
In our computations we refine such intervals until $|\gamma_2 - \gamma_1| < 10^{-100}$.

This approach has several advantages. 
First, it does not depend on the stability of the extrapolation procedure. 
Second, it is computationally more efficient, as it requires significantly fewer coefficients.

We use the following procedure to determine the $n^{\text{th}}$ root  of $S_1(\gamma)$. 
As an initial approximation ($i=0$), we take
$\gamma_0=\gamma_0(n)$
suggested by the asymptotic formula \eqref{Eq:gamma0asympt}. This approximation
works well for all $n\ge1$.
At each subsequent step, we employ the recurrence relation~\eqref{Eq:recursionak} 
and its derivative with respect to $\gamma$ 
to compute $a_k(\gamma_i)$ and $a'_k(\gamma_i)$ 
till we detect the onset of the monotone tail for $\hat a_k(\gamma_i)$ with
$k\ge k_{0}(\gamma_i)$.
 Then we perform one step of Newton's method:
\[
\gamma_{i+1} = \gamma_i - \frac{a_k(\gamma_i)}{a'_k(\gamma_i)}
\]
in order to find $\gamma$ with a smaller value of $a_k$.
This procedure is iterated until 
$|\gamma_{i+1} - \gamma_i| < 10^{-100}$, which is typically achieved in less
 than 10 steps.
We let $\gamma_*(n):=\gamma_{i+1}$.
As a final check, we verify that
\[
a_k\bigl(\gamma_*(n) + 10^{-101}\bigr) \, a_k\bigl(\gamma_*(n) - 10^{-101}\bigr) < 0
\]
for $k$ in the tail region.

\begin{figure}[htb]
    \centering
    \includegraphics[width=0.6\linewidth]{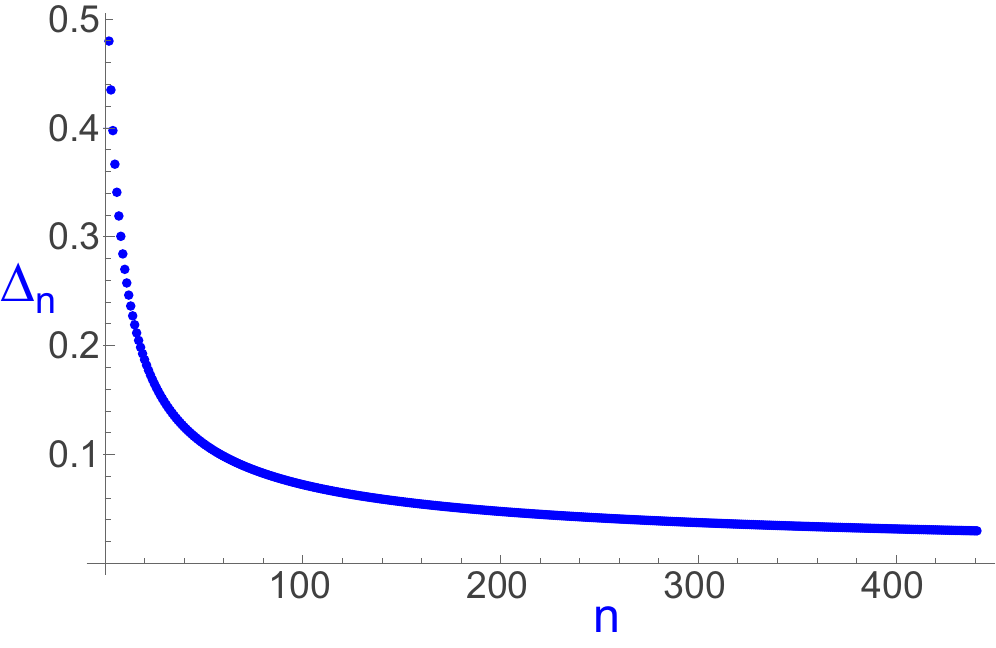}
    \caption{The first $440$ zeroes of $S_1$: deviation from the predicted value, $\Delta_n$,  in the units of the limit root spacing  vs $n$.}
    \label{Fig:zeroes}
\end{figure}

Figure~\ref{Fig:zeroes} shows the plot of the normalized deviation of the computed zeroes of $S_1(\gamma)$ from their
theoretical prediction,
\[\Delta_n=(\gamma_*(n)^{1/2} - \gamma_0(n)^{1/2})\,\sqrt{48}/\pi.
\]
The difference apparently  converge to zero as $n^{-1/2}$.
This conclusion confirms our theoretical prediction for the position of zeroes given by equation \eqref{Eq:gamma0asympt}.

\subsection{Method based on the variational equation}

A more direct computation of the Stokes constant is based on numerical evaluation of the difference 
between the stable and unstable solutions of the inner equation
on the negative imaginary axis and expressing this
difference in terms of the exponentially decaying solution
of the variational equation as defined in \eqref{Eq:S0}.

We construct an approximation for the unstable solution
$u^-$ to the inner equation using a method which is routinely used
in dynamical system for construction of unstable manifolds.
First we compute the coefficients of the formal solution~\eqref{Eq:uf} to the inner equation~\eqref{Eq:inner} using the recurrence relation~\eqref{Eq:recursionak}, and define the truncated series
\begin{equation}\label{Eq:un}
u_n(x)=\sum_{k=1}^{n} \frac{a_k}{x^{2k}}.
\end{equation}
Then we fix constants $x_R > \sigma_0 > 0$ and set the initial point at $x_0 = -x_R - \mathrm i \sigma_0$. 
Then we  solve numerically the Cauchy problem for the inner equation \eqref{Eq:inner} with initial conditions
\begin{equation}\label{Eq:initconditions}
u(x_0)=u_n(x_0), \quad u'(x_0)=u_n'(x_0), \quad u''(x_0)=u_n''(x_0), \quad u'''(x_0)=u_n'''(x_0),
\end{equation}
where $n$ is chosen according to the rule of truncation at the least term. 
The approximation to the unstable solution $u^-$ is numerically computed along the line $\Im(x) = -\sigma_0$ in the complex plane and provides a good approximation to $u^-(x)$
for $-x_R\le \Re x\le x_R$.

 The stable solution is defined using the reversible symmetry, $u^+(x)=
 \overline{ u^-(-\overline{x})}$. 
 On the imaginary axis the difference $\delta=u^+-u^-$ is purely imaginary:
 \[
 \delta(-\mathrm i\sigma_0) = u^-(-\mathrm i\sigma_0)-\overline{ u^-(\mathrm i\sigma_0)} =2 \mathrm i\Im u^-(-\mathrm i\sigma_0).
 \]
 Then we find coefficients of the formal solution 
 \eqref{Eq:deltaformal}
of the variational equation \eqref{EQ:innervareqn}
with $u=u_n$ and define
a partial sum
\begin{equation}\label{Eq:xim}
\xi_m(x)=\sum_{k=1}^{m} \frac{b_k \mathrm e^{-\mathrm i x}}{x^{k}}.
\end{equation}
Choosing $m$ from the optimal truncation rule for $x=-\rmi \sigma_0$
we obtain an approximation for the Stokes constant
\begin{equation}\label{Eq:S0ratio}
S_0\approx\frac{\delta(-\mathrm i \sigma_0 )}{
\xi_m(-\mathrm i \sigma_0)
}\;.
\end{equation}
This method yields results in excellent agreement with those obtained from the analysis of the coefficients. 
The theory predicts that the two definitions of the Stokes constant are related by
\[
S_0(\gamma) = -\pi \mathrm i\, S_1(\gamma).
\]
Our numerical experiments confirm this relation with an error below $10^{-14}$.

The method described in this section has a more direct and transparent connection to the Stokes phenomenon. 
However, achieving high precision is more challenging, as it requires the computation of highly accurate solutions to the Cauchy problem for the underlying ODE, which is computationally expensive.

The best precision achievable by this method is governed by the balance between two competing effects. 
On the one hand, the quantity $\delta$ decays exponentially when  $\sigma_0$ increases, while the accuracy of its evaluation is limited by the numerical precision of the solution to the Cauchy problem. 
On the other hand, the truncated-series approximation of the solution to the variational equation 
improves in accuracy as $\sigma_0$ increases and is reliable only for sufficiently large values of $\sigma_0$.


\subsection{Method based on the limit behaviour}

In this section, we describe an alternative method based on the definition of $S_0$ as the limit \eqref{Eq:S0def}, rather than on the ratio \eqref{Eq:S0ratio}. An important advantage of this method is that it does not require knowledge of the solution to the variational equation. However, it requires a very-high-precision Taylor-series integrator to reduce errors in solving the Cauchy problem in the complex time plane to an acceptable level.

\medskip

\textbf{Algorithm:}
\begin{enumerate}
\item For a given value of $\gamma$  we fix  $\sigma_0>0$ and
consider a  large  $x_\textrm{R}>0$ and define $x_0=-x_\textrm{R}-\mathrm i \sigma_0$. Then we numerically 
solve the Cauchy problem for the inner equation \eqref{Eq:inner} with initial conditions \eqref{Eq:initconditions}:
integrating the equation during time $t_\textrm{R}$,
we  obtain a very accurate approximation for $u^-(-\rmi \sigma_0)$ and its derivatives at this point.

\item
Then we choose $\sigma_{f}>\sigma_0$ and integrate the 
inner equation \eqref{Eq:inner} in the purely imaginary
direction substituting  $x=-\mathrm i \sigma$ and using
$\sigma$ as the new independent variable. 
We use the results of the previous step to generate
the initial condition for the Cauchy problem with $\sigma=\sigma_0$ and integrate the equation numerically
till $\sigma=\sigma_f$.

\item Once performed the numerical integration we use equation~\eqref{Eq:S0def} 
to get an approximation for the Stokes constant,
\begin{eqnarray}
S_0(\gamma) &=& \lim_{\Im x\to -\infty }\mathrm e^{\mathrm i x}\bigl(u^+(x)-u^-(x)\bigr)
\nonumber \\
\nonumber
&=& 2 \lim_{\sigma \rightarrow +\infty} \rme^{\sigma} \mathrm{Im}\left( u^+(-\rmi \sigma) - u^-(-\rmi \sigma) \right)
\\
&\approx& 2 \rme^{\sigma_f} \,\textrm{Im}(u^-(-\rmi \sigma_{f})).
\label{appstok}
\end{eqnarray}
\end{enumerate}

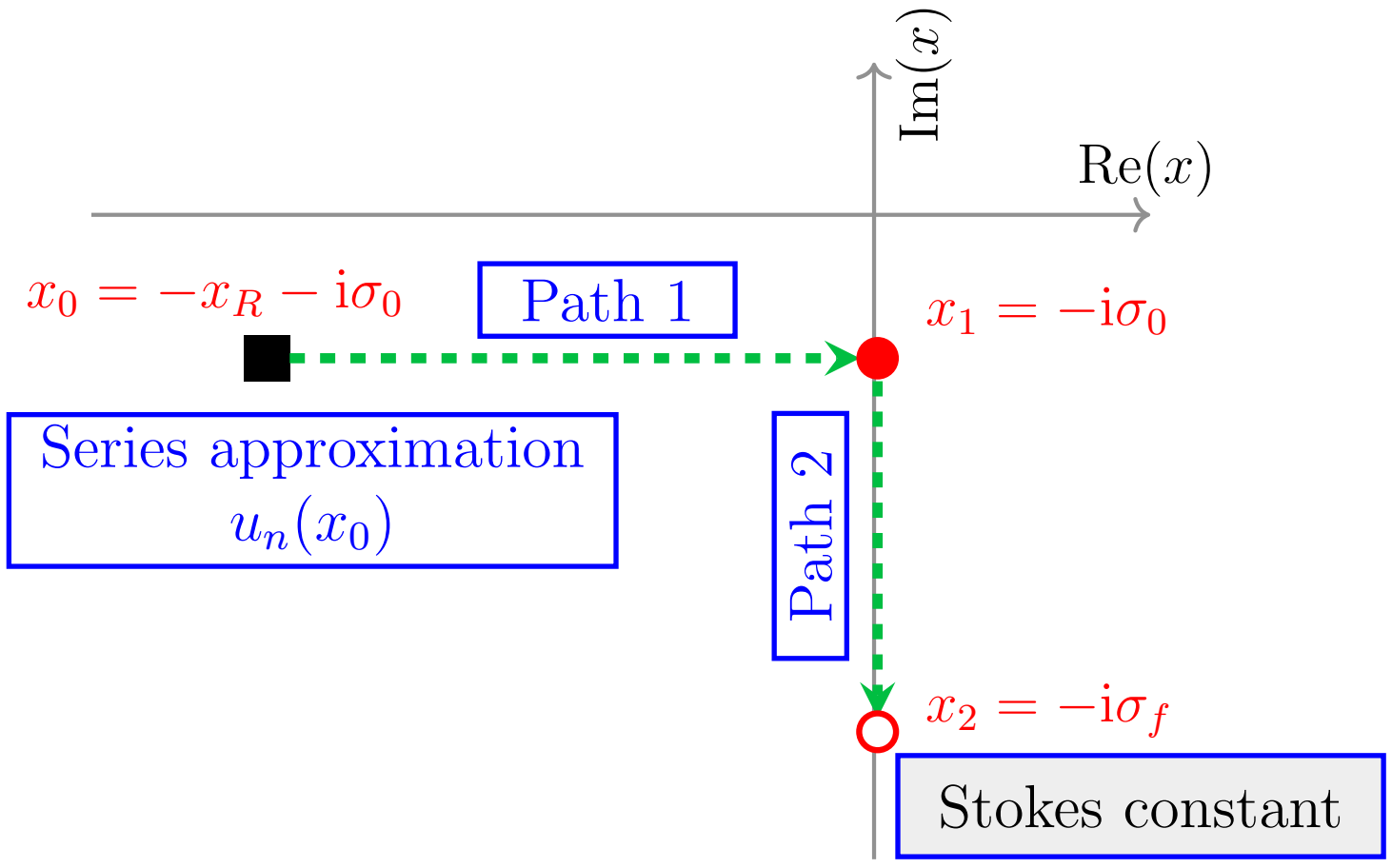
\begin{figure}[htb]
\centering
\begin{tikzpicture}[x=1pt,y=1pt,scale=0.7]
  \path[use as bounding box] (0,0) rectangle (348,225);

  \draw[axisgray,line width=0.8pt,->] (24,182) -- (323,182);>
  \draw[axisgray,line width=0.8pt,->] (245,0) -- (245,225);
  \node[anchor=west,scale=0.92] at (297,195) {$\operatorname{Re}(x)$};
  \node[rotate=90,anchor=west,scale=0.92] at (259,197) {$\operatorname{Im}(x)$};

  \fill[black] (67,135) rectangle (80,148);
  \draw[pathgreen,dashed,line width=2pt,->,>=stealth]
    (80,141.5) -- (241,141.5);

  \draw[pathgreen,dashed,line width=2pt,->,>=stealth]
    (246,135) -- (246,40);

  \fill[pointred] (246,141.5) circle (6pt);
  \draw[pointred,line width=1.2pt,fill=white] (246,36) circle (5.2pt);

  \node[
    draw=labelblue,
    text=labelblue,
    line width=1pt,
    inner xsep=8pt,
    inner ysep=3pt,
    anchor=south west
  ] at (133,147) {Path 1};

  \node[
    draw=labelblue,
    text=labelblue,
    line width=1pt,
    inner xsep=7pt,
    inner ysep=3pt,
    rotate=90,
    anchor=south west
  ] at (238,56) {Path 2};

  \node[
    draw=labelblue,
    text=labelblue,
    line width=1pt,
    align=center,
    minimum width=120pt,
    minimum height=30pt,
    inner sep=2pt,
    anchor=south west
  ] at (0,82) {Series approximation\\$u_n(x_0)$};

  \node[text=pointred,anchor=west,scale=0.94] at (0,160)
    {$x_0=-x_R-\mathrm i\sigma_0$};
  \node[text=pointred,anchor=west,scale=0.94] at (254,155)
    {$x_1=-\mathrm i\sigma_0$};
  \node[text=pointred,anchor=west,scale=0.94] at (254,42)
    {$x_2=-\mathrm i\sigma_{f}$};

  \node[
    draw=labelblue,
    fill=boxgray,
    text=black,
    line width=1pt,
    minimum width=96pt,
    minimum height=20pt,
    inner sep=2pt,
    anchor=south west
  ] at (251,0) {Stokes constant};
\end{tikzpicture}

\caption{\protect{\small
            Integration paths in the complex plane.
         \label{path}}}
\end{figure}

The integration path used in Steps 1 and 2 is sketched in Fig. \ref{path}.
\smallskip

The algorithm has been implemented
in multiple precision. 
The numerical tests have been done using an
optimized variable-step size variable-order Taylor series
integrator (the free software \textsf{TIDES}~\footnote{https://sourceforge.net/projects/tidesodes/}\cite{tides,BARRIO20117940}) in combination with the
GNU multiple precision libraries \textsf{GMP}~\cite{GMP} and
\textsf{MPFR}~\cite{MPFR}. Note that Taylor methods are the only ODE solvers able to reach very high precision levels \cite{BAILEY201210106,Barrio2015TS}. 

\begin{figure}[htb]
\centerline{\includegraphics[width=0.8\textwidth]{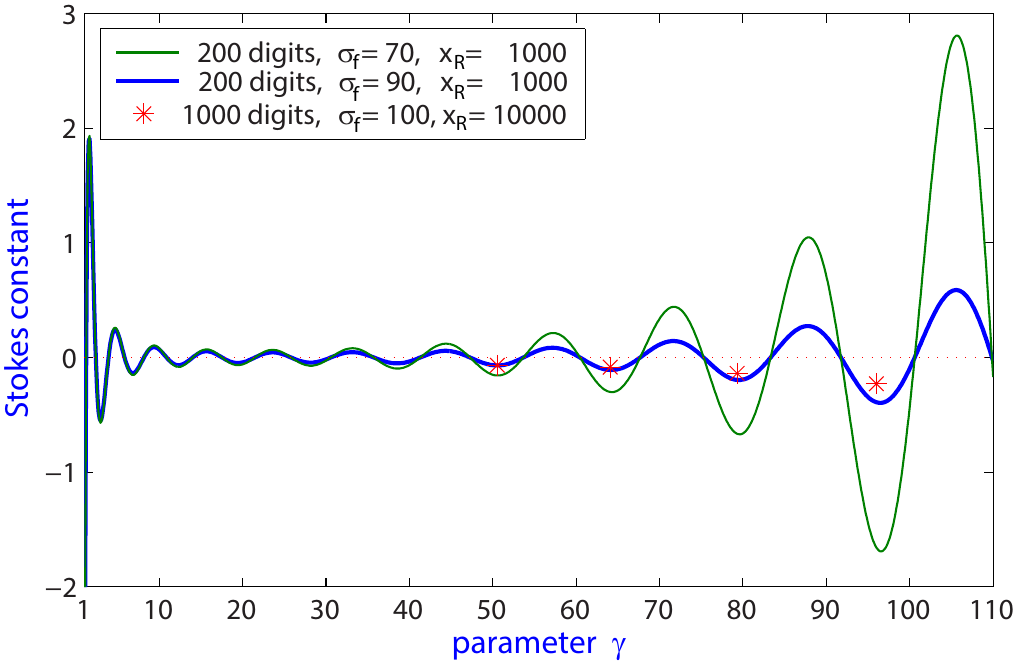}}
\caption{\protect{\small The estimates of Stokes constant $S_0(\gamma)$ for $\sigma_0 = 15$ using arithmetics with
         200 and 1000 decimal digits and different values of
         $\sigma_f$ and $x_\textrm{R}$.\label{tests1}}}
\end{figure}
\begin{figure}[htb]
\centerline{\includegraphics[width=1.\textwidth]{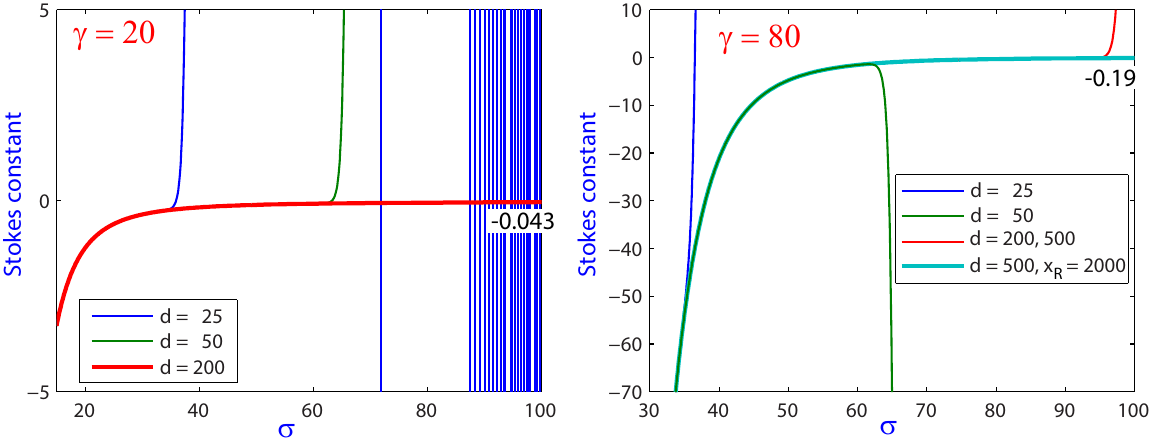}}
\caption{\protect{\small  Approximation of the
    Stokes constant approaching the limit value for $\gamma =
    20$ (left) and  $\gamma=80$ (right) using different precision levels (digits = $d$) and different
    initial times ($x_\textrm{R}$) for the numerical evaluation. \label{tests2}}}
\end{figure}

The above algorithm involves several parameters that must be chosen carefully: the displacement used to determine the initial condition
\[
x_0=-x_\textrm{R}-\rmi\sigma_0,
\]
the number of terms $n$ retained in the approximation 
of the initial conditions \eqref{Eq:un} (thus, the parameters associated with the first path are $x_\textrm{R}$, $\sigma_0$, and $n$), and the final imaginary time $\sigma_f$ along the second path. Maintaining high accuracy requires a careful choice of all these parameters and, consequently, extensive numerical testing.

In Fig.~\ref{tests1}, we show the estimates for the  Stokes constant computed with $\sigma_0=15$, which our tests indicate is sufficiently large, using either 200 or 1000 decimal digits of precision and several values of $\sigma_f$ and $x_\textrm{R}$. As $\gamma$ increases, obtaining an accurate   value of the Stokes constant becomes considerably more difficult. This is due both to the loss of accuracy in the numerical integration and/or the initial data, unless sufficiently high precision is used, and to the increasingly slow convergence, as illustrated in Fig.~\ref{tests2}.

For $\gamma=20$, computations with 25 decimal digits exhibit strong oscillations when the imaginary time $\sigma$ along the second path approaches $\sigma_f=36$, resulting in a highly inaccurate approximation of the Stokes constant. Increasing the precision to 50 digits delays the onset of the oscillations until $\sigma_f\simeq62$, whereas 200 digits are required to maintain adequate accuracy over the range considered. For $\gamma=80$, even higher precision is necessary, together with a sufficiently large value of the initial real displacement $x_\textrm{R}$.

High precision in the initial condition (IC) is crucial for obtaining an accurate value of the Stokes constant, although it is less important for locating its zeros $\gamma_{\textrm{num}}^{(i)}$, as shown in Fig.~\ref{tests1}. Table~\ref{table2} presents the number of correct digits in the numerical determination of the 21st zero of the Stokes constant, $\gamma_{\textrm{num}}^{(21)}$, for several values of $\sigma_f$, $x_\textrm{R}$, $n$, and the working precision. In each case, the quantity marked with an asterisk is the factor that limits the number of correct digits. These results show that an appropriate balance between the working precision and the accuracy of the initial condition is essential. This also requires a suitable choice of $n$ and $x_\textrm{R}$, since the initial point must be sufficiently far from the origin. Moreover, choosing a large value of $\sigma_f$ to assess convergence to the Stokes constant requires correspondingly higher precision; otherwise, numerical errors grow and may produce less accurate results than those obtained with a smaller value of $\sigma_f$.

\begin{table}[htb]
\begin{center}
\begin{tabular}{|ccccc|}
\hline
precision digits & $\sigma_f$ & $x_\textrm{R}$ & $J$ & correct digits
                   $\gamma_{\textrm{num}}^{(21)}$\\
\hline \hline &&&& \\[-0.3cm]
200 & 80${}^\ast$ & 1000 & 70 & 134 \\
200 & 80 & 1000 & 40${}^\ast$ & 122 \\
200 & 80 & 1000 & 20${}^\ast$ & 4 \\
200 & 80 & 100${}^\ast$  & 70${}^\ast$ & 9 \\
200 & 60 & 1000 & 70 & 173 \\
200${}^\ast$ & 80${}^\ast$ & 5000 & 70 & 134 \\
100${}^\ast$ & 60 & 1000 & 35${}^\ast$ & 53 \\
50${}^\ast$ & 60 & 1000 & 20${}^\ast$ & 3
\\
\hline
\end{tabular}
\caption{\protect{\small Number of correct digits in the numerical
    determination of the 21st zero of the Stokes constant
    ($\gamma_{\textrm{num}}^{(21)}$) for various values of
    $\sigma_f$, $x_\textrm{R}$, $J$ and precision digits. \label{table2}}}
\end{center}
\end{table}

\begin{table}[htb]
\hspace*{-2.8cm }
{\footnotesize
\begin{tabular}{|c|l|}
\hline
$i$ & $\gamma_{\textrm{num}}^{(i)}$ \\
\hline \hline & \\[-0.2cm]
1  & \texttt{1.25} \\
2  & \texttt{2.401049919639515487394205386216120634813759442477543869728648711149740720627580e+00}\\
3  & \texttt{3.930781055958246541778147654830619930170526589023534774773406052607795653141820e+00}\\
4  & \texttt{5.852175445170598283928616134085791767270340934872156470428868607330389222422947e+00}\\
5  & \texttt{8.171452010273213282253122677012053565621676960419593282409730035589907076181711e+00}\\
6  & \texttt{1.089222463289445398270765765983813970333283788702073647715940847811222347356477e+01}\\
7  & \texttt{1.401680066326203665412939722195407894879580483167280665124849956936106476459307e+01}\\
8  & \texttt{1.754674662406976296481029163323071683210209381245467456232394192894622165194123e+01}\\
9  & \texttt{2.148317536108186990216898283507374638983245347717053332798303126444890472248817e+01}\\
10 & \texttt{2.582690598464499793658543836930884706796068421901201419110303533400660368996520e+01}\\
11 & \texttt{3.057855900716822187495640835544144125961522482541362008239185246549622460852201e+01}\\
12 & \texttt{3.573861586739780457409126844937310873553227634093720444728233382699098112907917e+01}\\
13 & \texttt{4.130745772377472838094582156638522661775479887634910347718985055263149623896948e+01}\\
14 & \texttt{4.728539161014315214863834028085565871951097889866300564301068908172230572247931e+01}\\
15 & \texttt{5.367266858684872995682633469609052176254847669868734674493271657528979338213281e+01}\\
16 & \texttt{6.046949665426738133424415888966264863522612829910863019332850383006651406152488e+01}\\
17 & \texttt{6.767605014178120514390923103987227950622176380835836981450981438083384645476685e+01}\\
18 & \texttt{7.529247666606744212098373234109673095059304533634144776289061390477385198353880e+01}\\
19 & \texttt{8.331890237638685534133062359466256029977746734129818523264416732222084919042499e+01}\\
20 & \texttt{9.175543596915918363502484254063757642763319622100450702561676249884517762602421e+01}\\
21 & \texttt{1.006021718029178275020516133573017787027302673555997412932160876185711739569063e+02}\\
22 & \texttt{1.098591923453379924653981838458246921764827176552002828979475464931035514144293e+02}\\
\hline
\end{tabular}
}
\caption{The  first 22 zeros of the Stokes constant
    with 80 correct digits (the numerical process has been done
    with 200 precision digits, $n=70$, $x_\textrm{R}=1000$ and
    $\sigma_f=80$).\label{table1}}
\end{table}

\begin{figure}[htb]
\centerline{\includegraphics[width=0.6\textwidth]{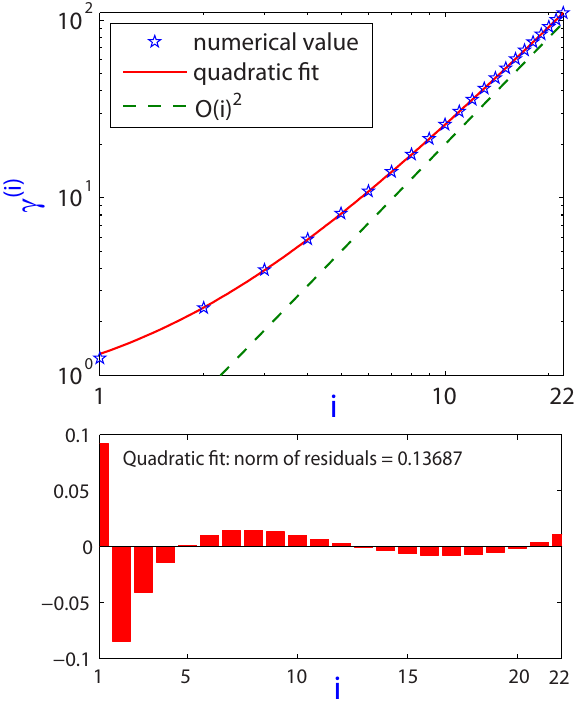}}
\caption{\protect{\small  First 22 zeros of the Stokes constant
    and a quadratic fit of the points in \texttt{log-log}
    scale.\label{roots2}}}
\end{figure}

Based on the preceding tests, and recognizing the need to balance computational cost against numerical accuracy, we selected the following parameter values: $d=200$ digits of precision, $n=70$, $x_\textrm{R}=1000$, and $\sigma_f=80$. Table~\ref{table1} lists the first 22 zeros of the Stokes constant, each computed to 80 correct digits. These values were independently verified using a higher-precision computation with $d=500$ digits of precision, $n=100$, $x_\textrm{R}=5000$, and $\sigma_f=90$.

The zeros were determined using Brent's method (see~\cite{Brent71}), which combines bisection, the secant method, and inverse quadratic interpolation. A quadratic least-squares fit to the computed zeros (Fig.~\ref{roots2}) yields
\[
 \gamma_{\textrm{num}}^{(i)} \approx 0.7391 + 0.4646\,i
 + 0.20431\,i^2, \qquad
 \|\textrm{error}\|_{2} = 0.13687.
\]
As shown in Fig.~\ref{roots2}, the quadratic approximation provides an excellent fit, with the error decreasing as $i$ increases.

Finally, the Stokes constant vanishes if and only if
\[
\Im u^-(-\rmi\sigma)=0 \qquad \text{for all } \sigma\geq\sigma_0.
\]
For such parameter values, $u^-$ is real along the imaginary axis. Therefore, if the aim is solely to locate the zeros, it suffices to integrate along the first path. The accuracy of a computed zero can then be verified by confirming that $\Im u^-(-\rmi\sigma_0)$ changes sign across the endpoints of a sufficiently small interval in~$\gamma$ containing the zero.

\section{Conclusions}

We have shown that the existence of a symmetric embedded soliton is determined by the zeros of a single parameter-dependent constant. This constant is the proportionality coefficient in the relation expressing the difference between the stable and unstable solutions in terms of a normalized solution of the variational equation. Using complex matching, we approximated this coefficient in terms of a Stokes constant associated with the inner problem. For large values of the parameter $\gamma$, the inner equation itself becomes singularly perturbed; we therefore applied the complex matching procedure a second time to derive an oscillatory asymptotic formula for the Stokes constant. We proved that this constant has countably many zeros that are asymptotically periodic in $\gamma^{1/2}$. Finally, we verified these asymptotic predictions by means of high-precision numerical computations.

\section{Acknowledgements}

The authors are grateful to Professor Alan Champneys, who introduced us to the problem of embedded solitons and brought together the team that completed the work presented in this paper.

RB has been supported by the project PID2024-156032NB-I0 funded by MICIU/AEI/10.13039/501100011033 and by ERDF, EU. JTL has been supported by the grant PID2024-158570NB-I00 funded by MICIU/AEI/ and ERDF, UE, 
and also by the Spanish State Research Agency, through the Severo Ochoa and María de Maeztu Program for Centers and Units of Excellence in R\&D (CEX2020-001084-M).

\bibliographystyle{plain}
\bibliography{biblio}

@article {Brent71,
    AUTHOR = {Brent, R.P.},
     TITLE = {An algorithm with guaranteed convergence for finding a
              zero of a function},
   JOURNAL = {Comput. J.},
  FJOURNAL = {The Computer Journal},
    VOLUME = {14},
      YEAR = {1971},
     PAGES = {422--425},
      ISSN = {0010-4620},
   MRCLASS = {65H05},
  MRNUMBER = {0339475},
MRREVIEWER = {W. M. Kincaid},
}

@InCollection{Barrio2015TS,
author="Barrio, R.",
editor="Engquist, Bj{\"o}rn",
title="Taylor Series Methods",
bookTitle="Encyclopedia of Applied and Computational Mathematics",
year="2015",
publisher="Springer Berlin Heidelberg",
address="Berlin, Heidelberg",
pages="1465--1470",
doi="10.1007/978-3-540-70529-1_153",
}

@article{BARRIO20117940,
title = "Breaking the limits: The {T}aylor series method",
journal = "Applied Mathematics and Computation",
volume = "217",
number = "20",
pages = "7940--7954",
year = "2011",
doi = "https://doi.org/10.1016/j.amc.2011.02.080",
author = "R. Barrio and M. Rodriguez and A. Abad and F. Blesa",
}

@article{BAILEY201210106,
title = "High-precision computation: Mathematical physics and dynamics",
journal = "Applied Mathematics and Computation",
volume = "218",
number = "20",
pages = "10106--10121",
year = "2012",
doi = "https://doi.org/10.1016/j.amc.2012.03.087",
author = "D.H. Bailey and R. Barrio and J.M. Borwein",
}

@article {ChampneysGroves97,
	AUTHOR = {Champneys, A.R. and Groves, M.D.},
	TITLE = {A global investigation of solitary-wave solutions to a
	two-parameter model for water waves},
	JOURNAL = {J. Fluid Mech.},
	FJOURNAL = {Journal of Fluid Mechanics},
	VOLUME = {342},
	YEAR = {1997},
	PAGES = {199--229},
	ISSN = {0022-1120},
	MRCLASS = {76B25 (34C37 35Q35 76B15)},
	MRNUMBER = {1460660},
	MRREVIEWER = {Boris Buffoni},
	DOI = {10.1017/S0022112097005193},
}

@article{ChampneysMalomedYang,
  author={Champneys, A.R. and Malomed, B.A. and Yang, J. and Kaup, D.J.},
  title={Embedded solitons: solitary waves in resonance
with the linear spectrum},
  journal={Physica D},
  year={2001}, 
  volume={152},
  pages={340-354} }

@article{Champmeys01PhysD,
	AUTHOR = {Champneys, A. R. and Malomed, B.A. and Yang, J. and Kaup, D.
	J.},
	TITLE = {Embedded solitons: solitary waves in resonance with the linear
	spectrum},
	JOURNAL = {Phys. D},
	FJOURNAL = {Physica D. Nonlinear Phenomena},
	VOLUME = {152/153},
	YEAR = {2001},
	PAGES = {340--354},
	ISSN = {0167-2789},
	MRCLASS = {35Q51 (35Q53 35Q60 81V80)},
	MRNUMBER = {1837916},
	DOI = {10.1016/S0167-2789(01)00178-6},
}

@article{YangEmbedded,
  author={Yang, J. and  Malomed, B.A. and Kaup, D.J.},
  title={Embedded solitons in second-harmonic-generating systems}, 
  journal={Phys. Rev. Lett.},
  volume={83},
  year={1999}, 
  pages={1958–1961} }

@article{IoossLombardi04,
  author = {Iooss, G. and Lombardi, E.},
  title = {Normal forms with exponentially small remainder: application
              to homoclinic connections for the reversible {$0\sp
              {2+}i\omega$} resonance},
  journal = {C. R. Math. Acad. Sci. Paris},
  fjournal = {Comptes Rendus Math\'ematique. Acad\'emie des Sciences.
              Paris},
  volume = 339,
  year = 2004,
  number = 12,
  pages = {831--838},
  issn = {1631-073X},
  mrclass = {37G05 (34C20 37C80)},
  mrnumber = {MR2111718 (2005h:37106)},
}

@book{Lombardi00,
  author = {Lombardi, E.},
  title = {Oscillatory integrals and phenomena beyond all algebraic
           orders},
  series = {Lecture Notes in Mathematics},
  volume = 1741,
  note = {With applications to homoclinic orbits in reversible
          systems},
  publisher = {Springer-Verlag},
  address = {Berlin},
  year = 2000,
  pages = {xvi+413},
  isbn = {3-540-67785-2},
  mrclass = {34C23 (34C37 34E10 37C29 37G25 37G40 76B25)},
  mrnumber = {MR1770093 (2002f:34081)},
  mrreviewer = {B. L. J. Braaksma},
}

@article{TovbisPelinovsky06,
  author = {Tovbis, A. and Pelinovsky, D.},
  title = {Exact conditions for existence of homoclinic orbits
           in the fifth-order {K}d{V} model},
  journal = {Nonlinearity},
  fjournal = {Nonlinearity},
  volume = 19,
  year = 2006,
  number = 10,
  pages = {2277--2312},
  issn = {0951-7715},
  coden = {NONLE5},
  mrclass = {37K99 (34C37 34M30 34M40 35Q53 37C29)},
  mrnumber = {MR2260265 (2007i:37143)},
  mrreviewer = {Vassilios M. Rothos},
}

@misc{GMP,
        author = {Torbj{\"{o}}rn et al., G.},
        title = "{GNU} Multiple Precision Arithmetic Library 4.1.2",
        month = {December},
        year = 2002,
}

@article{MPFR,
 author = {Fousse, L. and Hanrot, G. and Lef\`{e}vre,
           V. and P{\'e}lissier, P. and Zimmermann, P.},
 title = {MPFR: A Multiple-precision Binary Floating-point Library with
          Correct Rounding},
 journal = {ACM Trans. Math. Softw.},
 issue_date = {June 2007},
 volume = 33,
 number = 2,
 month = jun,
 year = 2007,
 issn = {0098-3500},
 articleno = 13,
 url = {http://doi.acm.org/10.1145/1236463.1236468},
 doi = {10.1145/1236463.1236468},
 acmid = 1236468,
 publisher = {ACM},
 address = {New York, NY, USA},
}

@article{tides,
 author = {Abad, A. and Barrio, R. and Blesa, F.
                    and Rodr\'{\i}guez, M.},
 title = {Algorithm 924: {TIDES}, a {T}aylor Series Integrator for
          Differential EquationS},
 journal = {ACM Trans. Math. Softw.},
 issue_date = {November 2012},
 volume = 39,
 number = 1,
 month = nov,
 year = 2012,
 issn = {0098-3500},
 pages = {5:1--5:28},
 articleno = 5,
 numpages = 28,
 url = {http://doi.acm.org/10.1145/2382585.2382590},
 doi = {10.1145/2382585.2382590},
 acmid = 2382590,
 publisher = {ACM},
 address = {New York, NY, USA},
}

@article {CraigGroves94,
	AUTHOR = {Craig, W. and  GrovesSegur,  M.D.},
	TITLE = {Hamiltonian long-wave approximations to the water-wave problem},
	JOURNAL = {Wave Motion},
	VOLUME = {19},
	YEAR = {1994},
	PAGES = {367--389},

}

@article {GelfreichLerman,
    AUTHOR = {Gelfreich, V. and Lerman, L.},
     TITLE = {Separatrix splitting at a {H}amiltonian {$0^2i\omega$}
              bifurcation},
   JOURNAL = {Regul. Chaotic Dyn.},
  FJOURNAL = {Regular and Chaotic Dynamics. International Scientific
              Journal},
    VOLUME = {19},
      YEAR = {2014},
    NUMBER = {6},
     PAGES = {635--655},
      ISSN = {1560-3547,1468-4845},
   MRCLASS = {37J20 (37J45 70K50 70K70)},
  MRNUMBER = {3284605},
MRREVIEWER = {Mikhail\ B.\ Sevryuk},
       DOI = {10.1134/S1560354714060033},
       URL = {https://doi.org/10.1134/S1560354714060033},
}

@book {Kuznetsov,
    AUTHOR = {Kuznetsov, Y.A.},
     TITLE = {Elements of applied bifurcation theory},
    SERIES = {Applied Mathematical Sciences},
    VOLUME = {112},
   EDITION = {Fourth},
 PUBLISHER = {Springer, Cham},
      YEAR = {2023},
     PAGES = {xxvi+703},
      ISBN = {978-3-031-22006-7; 978-3-031-22007-4},
   MRCLASS = {37Gxx (34C23)},
  MRNUMBER = {4592585},
       DOI = {10.1007/978-3-031-22007-4},
       URL = {https://doi.org/10.1007/978-3-031-22007-4},
}

@book {Meyer,
    AUTHOR = {Meyer, K.R. and Offin, D.C.},
     TITLE = {Introduction to {H}amiltonian dynamical systems and the
              {N}-body problem},
    SERIES = {Applied Mathematical Sciences},
    VOLUME = {90},
   EDITION = {Third},
 PUBLISHER = {Springer, Cham},
      YEAR = {2017},
     PAGES = {xiii + 384},
      ISBN = {978-3-319-53690-3; 978-3-319-53691-0},
   MRCLASS = {37J05 (37N05 70F10 70F15 70G45 70H05)},
  MRNUMBER = {3642697},
       DOI = {10.1007/978-3-319-53691-0},
       URL = {https://doi.org/10.1007/978-3-319-53691-0},
}

\appendix

\section{Hamiltonian structure\label{Se:hamiltonian}}

The equation \eqref{Eq:mainODE} has two properties important for our analysis.
First, the equation has a time-reversing symmetry, i.e., if 
a function $u(t)$ satisfies the equation, then $u(-t)$ is also
a solution of the same equation.
Second, the equation has a first integral and is Hamiltonian
with a non-standard symplectic form.

Let us describe the Hamiltonian structure in more details.
Denoting by $u_k$ the $k^{\mathrm{th}}$-derivative of $u$, we
rewrite the equation~\eqref{Eq:mainODE} in the form of a system of four first order equations:
\begin{equation}\label{Eq:hamsys}
\begin{split}
&u_0'=u_1,\quad u_1'=u_2,\quad u_2'=u_3,\\
&a_2 u_3'+a_1u_2+ a_0 u+b_0 u^2+b_1 ({u_1}^2+2u_0 u_2)=0.
\end{split}
\end{equation}
This system is hamiltonian with the Hamilton function 
\begin{equation}\label{Eq:hamiltonian}
	h=\frac{a_2}{2}\left(2u_3u_1-u_2^2\right)
	+\frac{a_1}{2}u_1^2
	+\frac{a_0}{2}u_0^2
	+\frac{b_0}{3}u_0^3
	+b_1u_0u_1^2
\end{equation}
and the non-standard symplectic form
\begin{equation}\label{Eq:symplecticform}
	\omega=a_2 (du_3\wedge du_0+du_1\wedge du_2)+(a_1+2b_1 u_0)du_1\wedge du_0.
\end{equation}
Indeed, let $X_h$ denote the vector field in the right hand side of the system \eqref{Eq:hamsys}. A direct substitution shows that  
 $X_h$  satisfies the identity $dh(\xi)=\omega(\xi,X_h)$ 
for all vectors $\xi=\sum_{k=1}^4\xi_k\partial_{u_k}$
where $\partial_{u_k}$ denotes the differentiation with respect to the
variable $u_k$. Recall that $du_k(\partial_{u_j})=\delta_{kj}$ where $\delta_{kj}$
is the Kronecker delta.

Soliton solutions belong to the zero level set of the first integral $h$.
Note that the equality $h=0$ can be interpreted
as an ordinary differential equation of order three. The equation has to be
satisfied by every solution of \eqref{Eq:mainODE} that converges to zero at infinity.

We note that the symplectic form can be easily transformed into the canonical form usually used in applications of the classical mechanics.

\section{Borel-Laplace transform and the Stokes constant\label{Se:resurgent}}

The choice of normalisation for the coefficients $a_k$ in the definition of $S_1$ as well as its connection to $S_0$ 
deserve an explanation. The factorial growth rate for
$a_k$ can be easily
proved with the help of the induction and 
recurrent  formulae \eqref{Eq:recursionak}.
The relation between $S_0$ and $S_1$  can be explained
provided with the help
from Ecalle's theory of resurgent functions. 
At the centre of this argument is the observation that the functions $u^+$ and
$u^-$ coincide with Borel-Laplace sums of the formal series $u_{\text{formal}}$.
In this interpretation the difference between $u^+$ and $u^-$ is due to
different choices of directions for Laplace integrals.

Let us discuss this point in more details.
It is possible to look for a solution of the equation~\eqref{Eq:inner}
in the form of a Laplace integral in a direction~$\alpha$, 
\begin{equation}\label{Eq:ualpha}
u_\alpha(x)=\int_0^{\mathrm{e}^{\mathrm{i}\alpha}\infty} e^{-sx}v(s)ds,
\end{equation}
where the integration path is a ray that forms the angle $\alpha$ with the positive real axis. We chose $\alpha=0$ for $u^+$
and $\alpha=\pi$ for $u^-$ and assume that $v$ is the same analytic functions in both cases.

Suppose that near the origin
\[
v(s)=\sum_{k=1}^{\infty}\frac{a_ks^{2k-1}}{(2k-1)!}.
\]
The growth rate of $a_k$ implies that 
 the series in the right hand side converge for $|s|<1$.
 Let $(u*v)(s)=\int_0^s u(t)v(s-t) dt$ be the convolution of functions.
 The Laplace transforms maps the convolution into the product and multiplication by $s$ into $\partial_x$. Therefore the Taylor series
 formally satisfy the integral equation 
\[
(s^4+s^2)v(s)+(v*v)(s)+\gamma \bigl((s v)*(s v)\bigr)(s)+2\gamma \bigl(v*(s^2 v)\bigr)(s)=0
\]
and, since the series converge, $v$ is an analytic solution of the same equation.

The integral equation can be used to prove that the analytic continuation of $v$
along any ray with $\alpha\ne \pm\frac\pi2\pmod {2\pi}$ does not grow faster
than an exponential function. We conclude that 
the Laplace integrals converge and $u_\alpha$ satisfies equation~\eqref{Eq:inner}. 
Then integrating by parts we can easily check that $u_\alpha\asymp u_{\text{formal}}$ as 
 the asymptotic behaviour of the Laplace integral for large $x$ is
determined by the Taylor series of $v$ around $s=0$.
The domain of validity of the asymptotic relation depends on $\alpha$,
and we conclude that our choices of $\alpha$ indeed lead to representation of the stable and unstable solution $u^\pm$
in the form of the Laplace integrals. 
Then  for all sufficiently large $x$ located in a neighbourhood
of the negative imaginary semi-axis we can write
\begin{equation}\label{Eq:deltaLaplace}
\delta (x)=u^+(x)-u^-(x)= u_\alpha(x)-u_{\pi-\alpha}(x)=
\int
_{-\mathrm{e}^{-\mathrm{i}\alpha}\infty}
^{ \mathrm{e}^{ \mathrm{i}\alpha}\infty} \mathrm e^{-sx}v(s)\,ds
\end{equation}
for some $\alpha>0$.

The upper bound on the growth rate of $a_k$ shows that
the series converge for $|s|<1$. In general the radius
of convergence is  equal to one. Indeed,
the integral equation can be rewritten in the form
\[
v=-\frac{v*v+\gamma (s v)*(s v)+2\gamma v*(s^2 v)}{s^2(s^2+1)}
\]
and consequently $v$ goes to infinity when $s\to\pm\rmi$ (unless
the integrals in the numerator vanish for some $\gamma$).
 Moreover, the dominant part of
the singularity is a simple pole.
Ecalle's theory suggests that the singularity should be of the form
\begin{equation}\label{Eq:vLaplace}
v(s)=\frac{C_0}{2\pi \mathrm i(s-\mathrm i)}+
\frac{\log(s-\mathrm i)}{2\pi\mathrm i}\sum_{k\ge0}c_k(s-\mathrm i)^k+\text{a regular function}(s-\mathrm i).
\end{equation}
Then deforming the path of \eqref{Eq:deltaLaplace}
into the upper half plane we conclude that
\[
\delta\asymp v_f=\mathrm e^{-\rmi  x}
\left(C_0+\sum_{k\ge 0}\frac{c_k k!}{x^{k+1}}\right).
\]
Comparing this expansion with \eqref{Eq:S0def} we conclude that
\[
S_0=C_0={2\pi \mathrm i}\lim_{s\to\lambda_0}(s-\mathrm i)v(s).
\]

\begin{prop} 
$S_0=\pi\mathrm i\,S_1$.
\end{prop}

A proof is based on the observation that the nearest singularity of a function 
is determined by the growth rate of coefficients of its Taylor series.
Since the function $v$ is real analytic, its poles are symmetric.
Moreover, we note that it is an odd function with real Taylor coefficients.
It follows that $v$ is purely imaginary on the imaginary axis.
Consequently, $\Re C_0=0$. Then
we note that
\[
\frac{C_0}{2\pi \mathrm i(s-\mathrm i)}+\frac{C_0}{2\pi \mathrm i(s+i)}
=
\frac{C_0 s}{\pi \mathrm i(s^2+1)}
=
\frac{C_0}{\pi \mathrm i}\sum_{k\ge 0}(-1)^{k}s^{2k+1}.
\]
Comparing this expansion with the Taylor series for $v$ and using 
 equation \eqref{Eq:vLaplace} we get
\[
\frac{C_0(-1)^k}{\pi \rmi}:\frac{a_k}{(2k-1)!}\to1.
\]
Consequently,
\[
S_1=\lim_{k\to\infty} \frac{(-1)^ka_k}{(2k-1)!}=\frac{C_0}{\pi\mathrm i}.
\]

\end{document}